\magnification=1200
\overfullrule=0mm

\font\tenrom=cmr10
\font\fiverom=cmr10 at 5pt
\font\sevenrom=cmr10 at 7pt
\font\eightrom=cmr10 at 8pt
\font\ninerom=cmr10 at 9pt
\newfam\romfam
\scriptscriptfont\romfam=\fiverom
\textfont\romfam=\tenrom
\scriptfont\romfam=\sevenrom

\font\nineromtt=cmtt10 at 9pt

\font\cyr=wncyr10
\font\ninecyr=wncyr10 at 9pt

\font\cyri=wncyi10

\newtoks\auteurcourant      \auteurcourant={\hfil}
\newtoks\titrecourant       \titrecourant={\hfil}

\newtoks\hautpagetitre      \hautpagetitre={\hfil}
\newtoks\baspagetitre       \baspagetitre={\hfil}

\newtoks\hautpagegauche   \newtoks\hautpagedroite 
  
\hautpagegauche={\ninecmr\rlap{\folio}\eightcmr\hfil\the\auteurcourant\hfil}
\hautpagedroite={\eightcmr\hfil\the\titrecourant\hfil\ninecmr\llap{\folio}}

\newtoks\baspagegauche      \baspagegauche={\hfil} 
\newtoks\baspagedroite      \baspagedroite={\hfil}

\newif\ifpagetitre          \pagetitretrue

\headline={\ifpagetitre\the\hautpagetitre
            \else\ifodd\pageno\the\hautpagedroite
             \else\the\hautpagegauche
              \fi\fi}

\footline={\ifpagetitre\the\baspagetitre\else
            \ifodd\pageno\the\baspagedroite
             \else\the\baspagegauche
              \fi\fi
               \global\pagetitrefalse}

\def\raggedbottom{\topskip 10pt plus 36pt\r@ggedbottomtrue}

\newcount\notenumber \notenumber=1
\def\note#1{\footnote{$^{{\the\notenumber}}$}{\eightrom {#1}}
\global\advance\notenumber by 1}

\auteurcourant={Philippe NUSS}
\titrecourant={GALOIS-AZUMAYA EXTENSIONS AND THE BRAUER-GALOIS GROUP}

\font\goth=eufm10
\font\tengoth=eufm10
\font\fivegoth=eufm10 at 5pt
\font\sevengoth=eufm10 at 7pt
\newfam\gothfam
\scriptscriptfont\gothfam=\fivegoth
\textfont\gothfam=\tengoth
\scriptfont\gothfam=\sevengoth
\def\goth{\fam\gothfam\tengoth}

\font\tengoth=eufm10
\font\fivegoth=eufm10 at 5pt
\font\sevengoth=eufm10 at 7pt
\newfam\gothfam
\scriptscriptfont\gothfam=\fivegoth
\textfont\gothfam=\tengoth
\scriptfont\gothfam=\sevengoth
\def\goth{\fam\gothfam\tengoth}
\font\ninecmr=cmr10 at 9pt
\font\eightcmr=cmr10 at 8pt
\font\ninecmsl=cmsl10 at 9pt
\font\ninecmbx=cmbx10 at 9pt
\font\eightcmbx=cmbx10 at 8pt
\font\eightcmsl=cmsl10 at 8pt

\font\twelvecmbx=cmbx10 at 12pt

\font\ninecmsl=cmsl10 at 9pt

\font\tensymb=msxm9
\font\fivesymb=msxm9 at 5pt
\font\sevensymb=msxm9 at 7pt
\newfam\symbfam
\scriptscriptfont\symbfam=\fivesymb
\textfont\symbfam=\tensymb
\scriptfont\symbfam=\sevensymb

\font\tensymbol=cmsy10
\font\fivesymbol=cmsy10 at 5pt
\font\sevensymbol=cmsy10 at 7pt
\newfam\symbolfam
\scriptscriptfont\symbolfam=\fivesymbol
\textfont\symbolfam=\tensymbol
\scriptfont\symbolfam=\sevensymbol
\def\symbol{\fam\symbolfam\tensymbol}

\newfam\bffam
\scriptscriptfont\bffam=\fivebf
\textfont\bffam=\tenbf
\scriptfont\bffam=\sevenbf
\def\bf{\fam\bffam\tenbf}

\font\tenmathb=cmmib10
\font\eightmath=cmmi10 at 8pt

\def\CB{\hbox{\bf C}}
\def\FB{\hbox{\bf F}}
\def\HB{\hbox{\bf H}}
\def\QB{\hbox{\bf Q}}
\def\RB{\hbox{\bf R}}
\def\RBIc{\hbox{\eightcmbx R}}
\def\ZB{\hbox{\bf Z}}
\def\ZBIc{\hbox{\eightcmbx Z}}

\def\ZC{\hbox{$\cal Z$}}

\def\H{\hbox{$\cal H$}}

\def\deltas {\hbox{$\partial$}}
\def\Del {\hbox{$\Delta$}}
\def\id{\hbox{\rm id}}
\def\Spec{\hbox{\rm Spec}}

\def\tr{\hbox{\rm tr}}

\def\dm{ \hfill {\fivesymb {\char 3}}}
\def\lr {\hbox{$\ \longrightarrow\ $}}
\def\ot {\hbox{$\otimes$}}
\def\pa{\S\kern.15em}
\def\Dem{\noindent {\sl Proof.}$\, \,$}

\def\limproj{\mathop{\oalign{lim\cr
\hidewidth$\longleftarrow$\hidewidth\cr}}}

\def\hfll#1#2{\smash{\mathop{\hbox to 10mm{\rightarrowfill}}
\limits^{\scriptstyle#1}_{\scriptstyle#2}}}

\def\vfl#1#2{\llap{$\scriptstyle #1$}\left\downarrow
\vbox to 3mm{}\right.\rlap{$\scriptstyle #2$}}

\def\diagram#1{\def\normalbaselines{\baselineskip=0pt
\lineskip=10pt\lineskiplimit=1pt}  \matrix{#1}}

\def\cross{\hbox{\symbol{\char 2}}}
\def\lcross{\hbox{\sevensymbol{\char 2}}}
\def\mmu{\hbox{\tenmathb{\char 22}}}

\def\Aut{{\rm Aut}}

\noindent {\twelvecmbx GALOIS-AZUMAYA EXTENSIONS AND
THE}
\medskip
\noindent {\twelvecmbx  BRAUER-GALOIS GROUP OF A COMMUTATIVE RING}

\vskip 25pt

\noindent {PHILIPPE NUSS}

\noindent {\ninerom Institut de Recherche Math\'ematique Avanc\'ee, Universit\'e Louis-Pasteur et CNRS, 
7, rue Ren\'e-Descartes,
67084 Strasbourg Cedex, France. e-mail: {\nineromtt nuss@math.u-strasbg.fr}}


\vskip 5pt
\noindent {\ninecmbx Abstract.}  {\ninecmr 
For any commutative ring {\ninecmsl R}, we introduce a group attached to {\ninecmsl R}, the 
{\ninecmsl Brauer-Galois group of {\ninecmsl R}}, defined to be 
the subgroup of the Brauer group of {\ninecmsl R}
consisting of the classes of the Azumaya {\ninecmsl R}-algebras which can be represented,
via Brauer equivalence, by a Galois extension of {\ninecmsl R}. We compute this group 
for some particular
commutative fields.

}

\vskip 5pt
\noindent {\ninecmbx Key-words:} {\ninecmr 
noncommutative ring, Galois-extension, Azumaya algebra, quaternion, Brauer group.}

\vskip 25pt
\noindent {\bf Introduction.}

\noindent 
Galois extensions of noncommutative 
rings were introduced in 1964 by Teruo Kanzaki [13].
These algebraic objects generalize
to noncommutative rings the classical Galois extensions of fields
and the Galois extensions of commutative rings due to Auslander and Goldman [1].
At the same time they also  turn out to be fundamental examples of 
Hopf-Galois extensions; these were first considered 
by Kreimer-Takeuchi [18] as a noncommutative analogue of the torsors
in algebraic geometry. Since Galois extensions are separable (Corollary 2.4) and since
the class of central Galois extensions $\psi : R \lr S$ over a fixed commutative
ground ring $R$ behaves well
under tensor product (Theorem 3.1),  
we may introduce a subgroup of the Brauer group of $R$, 
that we designate by {\sl Brauer-Galois group of $R$}. The purpose of the present paper
is to compute this object in some particular cases.
 
\smallskip
First of all, we need to recall the definitions and collect some basic properties of
Galois extensions. This is done in  section 1. 
In the second part, we investigate particular ring extensions such as
centralizing, Galois-Azumaya, 
Frobenius or separable extensions. 
We observe that Galois extensions
are Frobenius and
deduce from this fact that there exists a quantum object, called $\cal R$-matrix,
which is in a natural way attached  to any Galois-Azumaya extension $\psi$.
Therefore  $\psi$ gives rise, for any $n \geq 2$,
to a representation of the braid group $B_n$ (\pa \ 2.3).
The quaternion algebras furnish an archetypal family of examples of Galois-Azumaya extensions 
over their centre.
We show that the  ``$n$'th power norm residue algebras" --- a generalization of 
quaternion algebras originate from Class Field Theory --- are also Galois-Azumaya
extensions over their centre (Theorem 2.11). 
In the last section, we  prove that
the category of centralizing Galois extensions over a fixed ground ring $R$ 
 is stable  under tensor 
product (Theorem 3.1). This property enables us to introduce ${\rm Br}_{\rm Gal}(R)$, 
the Brauer-Galois group of a commutative ring $R$: it is the subgroup
of the Brauer group of $R$ 
consisting of those
classes $\xi \in {\rm Br}(R)$ for which there exists a finite group $G$ and 
a $G$-Galois-Azumaya extension $\psi : R \lr S$, such that $\xi$ is equal to the
class $[S]$ of $S$ modulo stable isomorphisms.
The computation of this object is in general at least as difficult 
as the determination of the classical Brauer group is.
However, when $n$ is a positive integer
and the ring $R$ is a field $F$ with  some 
arithmetic conditions with respect to $n$ satisfied, we may use results by 
Merkurjev and Suslin in order to identify  the $n$-torsion
of ${\rm Br}_{\rm Gal}(F)$ with that of ${\rm Br}(F)$
(Theorem 3.4). The main point lies in
the fact that the $n$'th power norm residue algebras
generate the $n$-torsion
of ${\rm Br}_{\rm Gal}(F)$ and are Galois-Azumaya
extensions over their centre. As a consequence, for some particular fields
of characteristic zero the Brauer group and the Brauer-Galois group coincide (Corollary 
3.5).

\bigskip
\noindent {\bf 1. Reminder on noncommutative Galois extensions.}
\smallskip
\noindent 
Let $k$ be a commutative ring, fixed throughout the paper
(possibly $k = \ZB$,  the ring of integers). By {\sl algebra} we mean an
associative  unital $k$-algebra.
A {\sl division algebra} is either a commutative field or a skew-field. 
If $R$ is an algebra,
an {\sl $R$-ring} $S$ is an algebra $S$ coming along with a morphism of algebras 
$\psi : R\lr S$. We call $\psi : R\lr S$ a {\sl (ring-)extension}. 
A morphism of extensions from $\psi : R\lr S$ to $\psi ' : R'\lr S'$
is a pair $(\varphi, \tilde \varphi)$,
where $\varphi: S \lr S'$ and $\tilde \varphi: R \lr R'$
are two morphisms of algebras  verifying
$\varphi \circ \psi  = \psi ' \circ \tilde \varphi $.
An $R$-ring $S$ carries always the $R$-bimodule structure
$r\cdot s\cdot r' = \psi (r)s\psi (r')$
(with $r,r'\in R$ and $s\in S$).

For any extension $\psi : R\lr S$, the group $S\ot _RS$
inherits an $S$-bimodule  (hence, via $\psi$, an $R$-bimodule) structure
by $s_1\cdot (s_2 \ot s_3)\cdot s_4 = s_1s_2 \ot s_3s_4$, with $s_1, s_2, s_3, s_4 \in S$.  A tensor
$\eta \in S\ot _RS$ is called a
{\sl Casimir element} in $S\ot _RS$ if it is
symmetric in the $S$-bimodule $S\ot _RS$, that is,
if it verifies, for all $s\in S$, the equality $\eta s = s\eta$.

An $R$-ring $S$ is endowed with an $R$-algebra structure
if and only if $R$ is commutative and the
morphism $\psi $ factors through the centre $\ZC (S)$ of $S$.
An extension  $\psi : R\lr S$
is {\sl centralizing} 
if the ring $R$ is embedded in the centre $\ZC (S)$ of $S$ (in this event
 $R$ is commutative and  $S$ is an $R$-algebra).
The extension  $\psi : R\lr S$
is {\sl central} if there is a bijection between the ring $R$ and the centre $\ZC (S)$
of $S$ (this means that $S$ is a faithful $R$-algebra and that we can identify
$R$ with $R\cdot 1$). An
extension $\psi : R \lr S$ is said to be {\sl commutative} if the rings $R$ and 
$S$ are commutative.

\medskip
\noindent 1.1. {\sl The definition of noncommutative Galois-extensions.}
\smallskip
\noindent
Before we recall the definition
of noncommutative Galois-exten\-sions, 
we need to introduce some notations. Let $S$ be an algebra and
$G$ be a finite group. 
Denote by $S(G)$ the algebra of all maps from $G$ 
to $S$, the sum and the product of two maps being calculated pointwise. Let
$\delta _g : G \lr S$ be the Dirac function at the point $g\in G$. It is
defined by $\delta _g(h) = \deltas _{g,h}1$, for any $h \in G$ ($\deltas _{g,h}$ stands for the 
Kronecker symbol of $g$ and $h$).  View $S(G)$ as the free left
$S$-module of rank $\vert G \vert$ (the order of the group $G$)  
with basis 
$\{ \delta _g, g \in G\}$, 
and with the product 
$$(s\delta _g)(t\delta _h) = st\deltas _{g,h}\delta _g,$$ 
for  $s,t \in S$, and $g,h \in G$. The family $(\delta _g)_{g\in G}$ is then  
a collection of orthogonal idempotents, the sum of which 
$\displaystyle\sum_{{g\in G}} \delta _g$ being
equal to the unit of the ring  $S(G)$. 

Suppose now that $G$ acts by ring automorphisms on $S$. Denote by
$g(s)$ the result of the action of $g \in G$ on an element $s \in S$. 
Endow $S(G)$ 
with the $S$-bimodule structure given, for $g\in G$ and $s_g, s, t \in S$, 
by the equality
$$s\cdot \sum_{g\in G}s_g\delta _g \cdot  t = \sum_{g\in G}ss_gg(t)\delta _g. \eqno (1)$$ 

\medskip
For any subalgebra $R$ of $S$ contained in the algebra $S^G$ of 
the invariant elements of $S$ under $G$,
define the map $\Gamma _{S/R} : S \ot _RS \lr S(G)$ by 
$$\Gamma _{S/R} (s \ot t) = \sum _{g\in G}sg(t)\delta _g = 
\sum _{g\in G}s\cdot \delta _g \cdot  t,$$
with  $s \ot t \in S \ot _RS$. This map $\Gamma _{S/R}$, sometimes denoted by $\Gamma _{\psi}$
or simply by $\Gamma$, is a
morphism of $S$-bimodules (hence of $R$-bimodules). Equipped with all these
notations, we are now able to state the definition of a Galois extension.

\medskip

\noindent {\bf  Definition 1.1.} A {\sl $G$-Galois extension} $\psi : R \lr S$
is a  morphism of algebras $\psi : R \lr S$, together with a finite group $G$
acting by ring automorphisms 
on $S$ and trivially on  $\psi (R) $ such that:

\item{--} the morphism $\psi$ is faithfully flat (that is $S$ is, via $\psi$, a faithfully flat left 
$R$-module);

\item{--} the morphism of $S$-bimodules $\Gamma _{S/R} : S\ot _RS \lr S(G)$, called 
in this case {\sl Galois morphism of $S/R$}, 
is an isomorphism.

\noindent The group $G$ is called a {\sl Galois group of $\psi$}.
A $G$-Galois extension  $\psi : R\lr S$
is {\sl centralizing}\note{Here the terminology differs with the one adopted
in [27], where we called such  extensions {\eightcmsl central}.} (respectively {\sl central}, 
respectively {\sl commutative}) if the underlying extension
is centralizing (respectively central, 
respectively commutative).
A $G$-Galois extension  $\psi : R\lr S$
is said to be {\sl strict} if the order of the 
group $G$ is invertible in the ring $S$. 
\medskip

\noindent 1.2. {\sl Preliminary results.}
\smallskip
\noindent
 Galois extensions verify 
many nice properties. We restate those
which are necessary 
in the rest of this paper.

\item{1)} If $\psi : R \lr S$ is a  $G$-Galois extension, 
the algebra of invariants $S^G$ is exactly
$\psi(R)$ and the morphism $\psi$ is injective [20]. Thus one may identify $\psi(R)$ with $R$
and treat $\psi$ as an inclusion. 

\item{2)} The abelian group  $S$ is, via $\psi$, also faithfully flat as a right 
$R$-module, and  $R$ is a direct summand as well of the right $R$-module $S$
as of the left $R$-module $S$ ([27], Proposition 1.2).

\item{3)} The left $R$-module  $S$ is, via $\psi$, a left 
  $R$-progenerator, that
is a finitely generated
projective left   $R$-module which is also a left 
  $R$-generator ([27], Lemme 1.4).

\item{3')} The same sentence as 3) but with ``left" replaced everywhere by ``right".

\medskip
The proof of these  statements call in a crucial way on a morphism of 
$R$-bimodules $S \lr R$, the {\sl trace map}, 
which appears as
some kind of $R$-linear section of the inclusion $\psi$.
Let $\psi : R \lr S$ be a  $G$-Galois extension. Since
$S^G$ is exactly $R$, the map from $S$ to $R$ denoted by $\tr _{S/R}$ or $\tr$ and
given by $$\tr (s) = \sum _{g\in G}g(s)$$ is well defined.
It is a morphism of $R$-bimodules, called the {\sl trace map}, and clearly
verifies $\tr \circ \psi = \vert   G \vert \cdot \id _R.$ 

\medskip

\noindent {\sl Galois extensions as Hopf-Galois extensions.} Let $G$ be a finite group and
let $\H = k^G $ be the Hopf algebra with $k$-basis
$\{ \delta _g\} _{g \in G}$, with multiplication $\cdot$
and comultiplication $\Del _{\cal H}$
defined by the formulae
$$\delta _g \cdot \delta _{g'} = \deltas _{g,g'}\delta _{g}  \quad {\rm and}
\quad \Del _{\cal H} (\delta _g)=
\sum _{ab = g}\delta _a \ot \delta _b$$
(the unit in $\H = k^G$ is the element 
${\displaystyle 1 = \sum _{g \in G}\delta _g}$ and the counit is the map
defined by $\epsilon (\delta _g) = 0$ if $g\not = e$ and
$\epsilon (\delta _e) = 1$).
A Galois extension $\psi : R \lr S$ with Galois group $G$ is then 
an $\H$-Hopf-Galois extension (we refer to [18] for the definition), where $S$ is an $\H$-comodule algebra with the coaction
$\Del _S: S \lr S \ot \H$ given by
${\displaystyle \Del _S(s) = \sum _{g\in G}g(s) \ot \delta _g.}$ 
\medskip

\noindent {\bf  Examples 1.2.} The prototypical examples of Galois extensions are:

\item{--} The classical finite Galois extensions of commutative fields.

\item{--} The Galois extensions of commutative rings, introduced by Auslander and Goldman [1],
studied 
by Chase, Harrison, A. Rosenberg and Sweedler ([5], [6], see also [17]).

\item{--} {\sl The algebra of diagonal matrices.} For any  division ring $D$ 
and any non-negative integer $n$, 
the diagonal map
$\psi  : D \lr D^n$ defines a $\ZB/n\ZB$-Galois extension ([27], \pa \  2.2).

\item{--} {\sl The trivial Galois extensions.} Given 
any (non)commutative ring $R$ and any finite group  $G$,
then $G$ is realizable over $R$, that is one may construct a
(non)commutative Galois extension $\psi : R\lr S$ having $G$ as Galois group. 
Do it in the following way:
set $S = R(G)$, embed $R$ into $R(G)$ via the diagonal map
${\displaystyle \psi (r) = \sum _{g\in G}r\delta _g}$ ($r\in R$), 
and let $G$ act  on $R(G)$ by the formula
$$g(r\delta _h) = r\delta _{hg^{-1}},$$
for $r\in R$ and $g,h \in G.$ Such a Galois extension is called {\sl trivial}.
It is straightforward to see that when one identifies
both $R(G) \ot _RR(G)$ and $R(G)(G)$ with $R(G\cross G)$, then the isomorphism
$\Gamma : R(G\cross G) \lr R(G\cross G)$ deduced from the Galois isomorphism
$\Gamma _\psi$
is given by  $$\Gamma (r\delta _{(g,h)}) = r\delta _{(g,g^{-1}h)},$$
for any $r\in R$ and $(g, h) \in G\cross G$.

\item{--} {\sl  The quaternion algebras over a field of characteristic 
different from $2$.} Let $F$ be a commutative field
of characteristic different from 2. Fix two elements
$a,b \in F^{\lcross}$.
Denote by ${\displaystyle ({{{a,b}\over F}})}$ 
the quaternion algebra over $F$. It is obtained by
dividing the free associative $F$-algebra on two generators $i$ and $j$
by the relations $i^2 = a$, $j^2 = b$ and $ij = -ji.$ 
It is well known (see [30] for example) that 
${\displaystyle ({{a,b}\over F})}$ 
is a central simple algebra of dimension 4 isomorphic either
to the algebra of matrices 
$M_2(F)$ or to a skew field, depending whether the 
Hilbert symbol $(a,b)_F$ is equal to $1$ or to $-1$
(when $F$ is the field  $\RB$ of the real numbers
and $a = b = -1$, then 
${\displaystyle ({{-1,-1}\over \RB})}$ 
is the skew field  
$\HB$ of Hamilton's quaternions). 
The Klein's {\sl Vierergruppe} $V = (\ZB /2\ZB)^2$ acts on ${\displaystyle ({{a,b}\over F})}$
by $\alpha (i)  = i$, $\alpha (j)  = -j$, $\beta (i)  = -i$, $\beta (j)  = j$, 
where $\alpha$ and $\beta$ are two generators of
$V$. 
In [26] (Lemma 4.3.2)  we proved that the extension 
${\displaystyle F \lr ({{a,b}\over F})}$ is $V$-Galois (see also [27], \pa \ 2.3).

\item{--} {\sl  A counter-example: The quaternion algebras over a field of characteristic 
 $2$.} Let $F$ be a commutative field
of characteristic  2. Fix two elements
$a,b \in F$. Denote by 
${\displaystyle {H_{a,b} =  [{{a,b}\over F})} }$ the quaternion algebra over $F$
(see [2]).
It is obtained by
dividing the free associative $F$-algebra on two generators $e_1$ and $e_2$
by the relations 
$e_1^2 = e_1 + a$, $e_2^2 = b$ and $e_2e_1 = e_1e_2 + e_2$.
The ring $H_{a,b}$ is a skew field if and only if the polynomial
$X_0^2 + X_0X_1 +aX_1^2 + b(X_2^2 + X_1X_2 + aX_3^2)
\in F[X_0, X_1, X_2, X_3]$
has the only one root, namely $(0,0,0,0)$.
When $H_{a,b}$ is a skew field, with
$b = c^2$ a square in $F$, we have shown ([26], \pa \ 2.4) that the extension
$\psi : F \lr H_{a,b}$ can never be Galois.

\medskip

\noindent For the sake of completeness, we state again an
important result, proved in [27] (Th\'eor\`eme 1.9), which we need in the sequel: 

\medskip

\noindent {\bf  Theorem 1.3.} {\sl  Let $G$ be a finite group and  $\psi : R\lr S$
be a strict $G$-Galois extension. For any subgroup 
$H$ of $G$, denote by $U = S^H$ the ring of fixed elements of $S$ under $H$ 
and by $\theta : U\lr S$ the canonical inclusion map. The morphism $\theta$ is then a 
strict $H$-Galois extension. Moreover, if $H$ is normal in $G$, 
the canonical inclusion map
$\theta ' : R\lr U$ is  a strict 
$G/H$-Galois extension.}

\medskip
\noindent 1.3. {\sl The Galois element.}
\smallskip
\noindent
 Let  $\psi : R \lr S$ be a $G$-Galois extension. For
any  $g\in G$, denote by 
$\eta_g$ the element $\Gamma^{-1}(\delta_g)$ in
$S\ot _RS$. The identity (1) implies the equality   $$\eta _gs = g(s)\eta _g,$$
which holds in $S\ot _RS,$ for all $s\in S$. 
Denote by $e$ the neutral element
of the group $G$.
Then $\eta_e \in S\ot _RS$ is called the {\sl Galois element for $\psi $}; 
it is in particular a Casimir element.
 Fix  once and for all a decomposition $\displaystyle{\sum_{i=1}^m x_i\ot y_i}$ of the
Galois element $\eta_e \in S\ot _RS$.
The $2m$-tuple $(x_1, \ldots ,x_m$; $y_1, \ldots , y_m)$
is called a
{\sl Galois basis of $\psi $.}
The equality  $\Gamma (\eta_e) = \delta_e$
can be expressed by  ${\displaystyle\sum _{i=1}^mx_ig(y_i) =  \deltas _{g,e}.}$
We state now, without proof, the 
properties verified by the elements $\eta_g \in S\ot _RS$,
which we need here   (for details, see [27]).

The element $\eta _g$ can be recovered
from the Galois element
$\displaystyle{\eta _e =  \sum _{i=1}^mx_i\ot y_i}$, 
namely $\displaystyle{\eta _g =  \sum _{i=1}^mx_i\ot g^{-1}(y_i)}$,
which implies the equality 
$\displaystyle{\sum _{i=1}^mx_igh^{-1}(y_i) =  \deltas _{g,h}}$.
Finally the collection of all $\eta _g$ (for $g \in G$) form a ``partition of the unity" in the 
following sense: one has $\displaystyle{\sum _{g\in G}\eta _g = 1\ot 1.}$

\medskip

\noindent {\bf  Lemma 1.4.} {\sl Let $G$ be a finite group and $\psi : R \lr S$ 
be a $G$-Galois extension. 
The opposite morphism $\psi ^{\rm o} : R^{\rm o} \lr S^{\rm o}$ defines 
then a $G$-Galois extension,
 $G$ acting on $S^{\rm o}$ via $$g(s^{\rm o}) = (g(s))^{\rm o}.$$}

\Dem It is clear that $(S^{\rm o})^G = (S^G)^{\rm o} = R^{\rm o}$. By
[27] (Th\'eor\`eme 1.5), it remains  to show that the morphism 
$\Gamma _{\psi ^{\rm o}}$ from $ S^{\rm o} \ot _{R^{\rm o}}S^{\rm o}$ to $S^{\rm o}(G)$ given by
$$\Gamma _{\psi ^{\rm o}} (s^{\rm o} \ot t^{\rm o}) = \sum _{g\in G} s^{\rm o}g(t^{\rm o})\delta _g = 
 \sum _{g\in G} (g(t)s)^{\rm o}\delta _g,$$ is surjective.
To this end, choose a Galois basis   $(x_1, \ldots ,x_m$; $y_1, \ldots , y_m)$ of $\psi $. Then
$$\Gamma _{\psi ^{\rm o}} \bigl(\sum _{i = 1}^m y_i^{\rm o} \ot g^{-1}x_i^{\rm o}\bigr) = 
\sum _{g\in G}\bigl(\sum _{i = 1}^mhg^{-1}(x_i)y_i\bigr)^{\rm o}\delta _h,$$
which is equal to $\delta _g$, according to the identity (5') in [27].
\dm

\medskip
Notice that the quaternion algebras ${\displaystyle ({{a,b}\over F})}$
over  a commutative field $F$
of characteristic different from 2 are isomorphic to their opposite algebras.

\medskip
\noindent 1.4. {\sl  The categories ${\goth Gal}$ and ${\goth Gal}(G)$.}
\smallskip
\noindent
First introduce the category
${\goth Gal}$ of Galois extensions: its objects are 
the couples $(G, \psi : R \lr S)$, with
$G$ a finite group and $\psi : R \lr S$ a $G$-Galois extension;
a morphism in ${\goth Gal}$ from $(G, \psi : R \lr S)$ to  $(G', \psi ' : R' \lr S')$
is a triple $(f, \varphi, \tilde \varphi)$,
where $f : G'\lr G$ is a morphism of groups and 
($\varphi : S \lr S', \tilde \varphi : R \lr R'$)
is a morphism of extensions (that is 
$\varphi  \circ \psi  = \psi ' \circ \tilde \varphi$) verifying
$\varphi  \circ f(g') = g' \circ \varphi $, for all $g'\in G'$.
Necessarily $\tilde \varphi $ is then the restriction of $\varphi $ to $R$.
In a similar manner, define the full subcategories 
${\goth Galstr}$, ${\goth Galcent}$, ${\goth GalAzum}$
or ${\goth Galcom}$
of ${\goth Gal}$; their objects are respectively the
strict, centralizing, central or commutative Galois extensions. 

Let ${\goth Gpf}$ be the category of finite
groups and $\pi ^{gr}: {\goth Gal} \lr {\goth Gpf}$ be the contravariant 
functor defined
by $\pi ^{gr}(G, \psi : R \lr S) = G$. The fibre category of $\pi ^{gr}$ over
a fixed group $G$ is denoted by
${\goth Gal}(G)$; its objects are therefore the  $G$-Galois extensions
$\psi : R \lr S$; a morphism from $(\psi : R \lr S)$ to $(\psi ' : R' \lr S')$
in ${\goth Gal}(G)$ is a couple $(\varphi  : S \lr S', \tilde \varphi : R \lr R')$
of morphisms of algebras such that $\varphi $  is $G$-equivariant (that is
$\varphi  \circ g = g\circ \varphi $ for all $g \in G$) and
$\psi ' \circ \tilde \varphi  = \varphi  \circ \psi $.  Observe that 
$\tilde  \varphi$ is induced by restriction by 
$\varphi$. Starting with  ${\goth Gal}(G)$, one may define in an obvious way the full subcategories 
${\goth Galstr}(G)$, ${\goth Galcent}(G)$, ${\goth GalAzum}(G)$
or ${\goth Galcom}(G)$
of ${\goth Gal}(G)$.

\medskip

\noindent {\bf  Proposition 1.5.} {\sl Let $(\varphi , \tilde \varphi )$ 
be a morphism from $(\psi : R \lr S)$ to $(\psi ' : R' \lr S')$
in ${\goth Gal}(G)$. 
Then $$\tilde \varphi  \circ \Gamma = \Gamma '\circ (\varphi  \ot \varphi  ).$$
Here $\Gamma$ (respectively $\Gamma '$) stands for the Galois isomorphism
of $S/R$ (respectively of $S'/R'$), and $\tilde \varphi $ is
the morphism of left $R$-modules from  $S(G)$ to $S'(G)$ defined by 
$\displaystyle{\tilde \varphi  (\sum _{g\in G}s_g\delta _g) =}$
$\displaystyle{\sum _{g\in G}\varphi  (s_g)\delta _g.}$
}

\medskip
\noindent {\sl Proof.} For any element $s\ot t \in S\ot _RS$, one has the equality
$$\bigl(\Gamma '\circ (\varphi  \ot \varphi )\bigr)(s\ot t) \ =
\displaystyle{\sum_{g\in G}\varphi  (s)g\bigl(\varphi  (t)\bigr)\delta _g},$$ 
whereas $(\tilde \varphi  \circ \Gamma )(s\ot t)$
=
$\displaystyle{\sum_{g\in G}\varphi  (s)\varphi  \bigl(g(t)\bigr)\delta _g}$.
The proposition is then a consequence of the $G$-equivariance of $\varphi $. \dm

\medskip

\noindent {\bf  Corollary 1.6.} {\sl Let $(\varphi , \tilde \varphi )$ 
be a  morphism from $(\psi : R \lr S)$ to $(\psi ' : R' \lr S')$
in ${\goth Gal}(G)$. If $\eta _e$ (respectively $\eta'_e$)
is the  Galois element for $\psi $ (respectively for $\psi '$), 
then $$(\varphi  \ot \varphi ) (\eta _e) = \eta '_e.$$}

\noindent {\sl Proof.} One has $\bigl(\Gamma '\circ (\varphi \ot \varphi  )\bigr) (\eta _e) $
= $(\tilde \varphi   \circ \Gamma ) (\eta _e) $ = $\tilde \varphi   (\delta _e) $ = $ \delta _e$.
The result follows from the fact that $\Gamma '$ is an isomorphism. \dm

\bigskip
\noindent {\bf 2. Galois-Azumaya extensions.}
\smallskip
\noindent 2.1. {\sl Separable extensions.} 
\smallskip
\noindent
Let  $\psi : R \lr S$ be a morphism of algebras and
$S^{e}  = S\ot _k S^{{\rm o}}$  the enveloping algebra of $S$ over the ground ring $k$.
Define the multiplication $\bar \mu : S^e \lr S$ to be the $S^e$-linear map given by
$\bar \mu (s\ot t^{\rm o}) = st$, and denote by $\Omega _{\bar \mu} $ the kernel of $\bar \mu $. One easily
sees that $\Omega _{\bar \mu} $ is generated by the elements 
$s\ot t^{\rm o} - st\ot 1^{\rm o}$ indistinctly as a left or as a right $S$-module.  
An element $\eta \in S\ot _RS$ is then Casimir if and only if
$\Omega _{\bar \mu} \cdot \eta  = 0.$

\medskip

\noindent {\bf  Proposition 2.1.} {\sl Let $\psi : R \lr S$ be a morphism of algebras.
The following conditions are equivalent:

{\rm (i)} The multiplication $\mu :  S\ot _RS \lr S$
splits as an $S^{e}$-module morphism.

{\rm (ii)} There exists a Casimir element 
$\eta $ in $S\ot _RS$ such that $\mu (\eta ) = 1.$ }

\noindent 
{\sl If moreover $R$ coincides with the (commutative) ground ring $k$, 
one of the two previous conditions is equivalent  to the following assertion:

{\rm (iii)} The $k$-module $S$ is projective as a left $S^{e}$-module.}

\medskip

\noindent {\bf  Definition 2.2.} 
A morphism of algebras $\psi : R \lr S$ is called {\sl separable} if it satisfies the equivalent
conditions of Proposition 2.1.
A Casimir element 
$\eta $ in $S\ot _RS$ verifying $\mu (\eta ) = 1$ is called a {\sl separability element of $S.$} 
If moreover the extension $\psi$ is centralizing, then $S$ is a {\sl separable $R$-algebra.} 
Furthermore, if $\psi$ is central, then $S$ 
is a central separable $R$-algebra, or an {\sl Azumaya $R$-algebra}\note{We follow here
the terminology
proposed by Bourbaki [3]. Some authors deal only with Azumaya algebras over a field.
These are central simple algebras [7] (as an example take quaternion algebras over a field). Recall
that an Azumaya {\eightmath R}-algebra is simple if and only if {\eightmath R} is a field.}.

\medskip
\noindent {\sl Remarks:}

\noindent 1.-- One may also characterize separable extensions using derivations or
Hochschild cohomology (see for example [4], 1.3, Theorem 3]).

\noindent 2.-- When  $R$ is the commutative ground ring $k$, the separability element $\eta$ viewed as an
element of $S^{e}$ is necessarily an idempotent of the algebra
 $S^{e}$. This means: fix a decomposition ${\displaystyle \sum_{i=1}^mu_i \ot v_i}$
of $\eta \in S\ot _RS$ and set $ {\displaystyle e(\eta)  = \sum_{i=1}^mu_i \ot v_i^{\rm o} \in S^{e}}$.
Then $e(\eta)^2 = e(\eta)$. 
Indeed, $e(\eta)^2-e(\eta) = e\bigl(e(\eta)\cdot e - (1\ot 1^{\rm o})\eta\bigr) \in e(\Omega _{\bar \mu} \cdot \eta) = 0.$

\medskip
\noindent {\sl Proof of  Proposition 2.1.}

\noindent {\rm (i) $\Longrightarrow$ (ii)} 
Let $\sigma$ be an $S^{e}$-linear section  of $\mu$ and $\eta = \sigma (1)$. 
Then $\mu (\eta) = \mu \circ \sigma (1) = 1.$
For  $s \ot t^{\rm o} -st \ot 1^{\rm o} \in \Omega _{\bar \mu}$, one has
$(s \ot t^{\rm o} -st \ot 1^{\rm o})\cdot \eta $ = 
$(s \ot t^{\rm o})\sigma (1) - (st \ot 1^{\rm o})\sigma (1) $
 = $\sigma (st) - \sigma (st) = 0,$
hence $\Omega _{\bar \mu}\cdot \eta = 0.$

\noindent {\rm (ii) $\Longrightarrow$ (i)} 
Since $\Omega _{\bar \mu}\cdot \eta = 0$, one defines a map  $\sigma$ from $S$ 
to $S\ot _RS$ by
$\sigma (s) = (s \ot 1^{\rm o})\cdot \eta$. Thus  $\sigma (s) = s\cdot \eta = 
(1 \ot s^{\rm o})\cdot \eta = \eta\cdot s.$
Then $\mu  \circ \sigma (s) = \mu (s\cdot \eta) = s\mu (\eta) = s$,
and $\sigma$ is $S^{e}$-linear since
$\sigma ((s \ot t^{\rm o})\cdot u) $ = $\sigma (sut)$
= $(sut \ot 1^{\rm o})\cdot \eta $
= $(su \ot 1^{\rm o})(t \ot 1^{\rm o})\cdot \eta $
= $(su \ot 1^{\rm o})(1 \ot t^{\rm o})\cdot \eta $
= $(su \ot t^{\rm o})\cdot \eta $
= $(s \ot t^{\rm o})(u \ot 1^{\rm o})\cdot \eta $
= $(s \ot t^{\rm o})\sigma (u).$

\noindent {\rm (i) $\Longleftrightarrow$ (iii)} is  well known (see [7]).
\dm

\medskip

\noindent Examples of separable extensions are ([7], [15], [17]): 

\item{i)} A finite product $\prod L_i/K $ of separable finite commutative field extensions 
$L_i/K$  (recall that a field extension $L/K$ is separable if and only if the minimal polynomials of the elements of
$L$ have simple zeroes);

\item{ii)} The $n$-fold product $R \cross R  \ldots \cross R$ of a commutative ring $R$ 
is separable over $R$;

\item{iii)} Let $A$ be a commutative algebra and $\Sigma$ be a multiplicative subset
of $A$. The  localized ring $\Sigma ^{-1}A$ is separable over $A$;

\item{iv)} The ring $M_n(R)$ of all $n \cross n$-matrices over  a commutative ring $R$ 
is separable over $R$;

\item{v)} Let $G$ be a finite group whose order $n$ is 
invertible in a commutative ring $R$. The group 
algebra $R[G]$ is separable over $R$;

\item{vi)} Any Azumaya algebra.

\medskip

\noindent {\bf  Definition 2.3.} A morphism of rings
$\psi : R \lr S$ is called {\sl a
$G$-Galois-Azumaya extension} if $\psi$ is a $G$-Galois extension, $R$ is a commutative ring
and  $S$ is an
$R$-Azumaya algebra, that is a 
central separable $R$-algebra.
\smallskip

\noindent For instance, the quaternion extension 
${\displaystyle F \lr ({{a,b}\over F})}$ (Example 1.2)
is $V$-Galois-Azumaya, with $V = (\ZB /2\ZB)^2.$ 

\medskip
Galois-Azumaya extensions verify nice properties in view of cohomology theories.
On the one side, 
the group of automorphisms of an Azumaya algebra is well controlled ([17], IV).
For instance, if $S$ is an Azumaya algebra over a ring $R$, 
then $\Aut (S/R)$ fits into an exact sequence of groups
$1\lr R^{\lcross}\lr S^{\lcross} {\lr } 
\Aut (S/R) \lr {\rm Pic}(R)$, known as the Rosenberg-Zelinsky exact sequence [17]. Therefore,
if the Picard group ${\rm Pic}(R)$ is trivial (what happens
when $R$ is a local ring or a principal ideal domain for example),
one obtains the Skolem-Noether theorem, which asserts that any
automorphism of $S$ leaving $R$ pointwise invariant is inner. 
On the other side, any $G$-Galois-Azumaya (or, more generally,
centralizing $G$-Galois) extension $\psi : R \lr S$ provides a natural
crossed module of $G$-groups obtained from $S^{\lcross} {\lr } 
\Aut (S/R)$ which can be taken as coefficients 
for nonabelian hypercohomology [27]. 

\medskip
We restate now a well-known and straightforward result [20] that furnishes further examples 
of separable extensions.

\medskip
\noindent {\bf  Corollary 2.4.} {\sl Any Galois extension is
separable.}
\medskip

\noindent Indeed, if  $\psi : R \lr S$ is a $G$-Galois extension, the Galois element $\eta _e$ is a
separability element for $S$, hence $\psi$ is separable. 
{\sl Therefore a central Galois extension amounts to a Galois-Azumaya  extension.}
We remark now that Galois extensions own another interesting structure, 
that of a Frobenius extension.
\medskip

\noindent 2.2. {\sl Frobenius extensions.} 
\smallskip
\noindent
A morphism of algebras $\psi : R\lr S$ is called a {\sl Frobenius extension} 
if there exists a morphism $\tau : S \lr R$ of $R$-bimodules
and a finite set of couples $(u_i, v_i)_{i=1, \ldots, m}$, with $u_i, v_i \in S$ verifying
for all $s\in S$ the {\sl Frobenius normalizing conditions}, that is
$$\sum _{i=1}^m u_i\tau (v_is) = s \ \ \ \ {\hbox {\rm and}}\ \ \ \ \sum _{i=1}^m \tau (su_i)v_i = s.$$

\noindent The collection $((u_i, v_i)_{i=1, \ldots, m})$ 
 allows us to construct
the tensor $\displaystyle{ \eta =  \sum _{i=1}^mu_i\ot v_i \in S\ot _RS}$,
called the {\sl Frobenius element of the extension $\psi : R\lr S$}.
The datum $((u_i, v_i)_{i=1, \ldots, m}, \tau)$ or $(\eta, \tau)$
is a {\sl Frobenius system}.
If $R$ lies in the centre of $S$, the extension
$\psi : R\lr S$ is {\sl centralizing Frobenius}
($S$ is then a  {\sl Frobenius $R$-algebra}).

\medskip

As examples of Frobenius extensions, let us mention, with Kadison, morphisms of integral group
algebras $\psi : \ZB [H] \lr \ZB [G]$ (where $H$ is a subgroup of $G$ of finite index
$\ell$, such that $G = \coprod g_iH$, with $g_1 = e, g_2, \ldots , g_{\ell} \in G$) ([12], Example 1.7),
and the  Frobenius extensions
obtained from certain types of von Neumann $*$-algebras ([12], Example 1.8).
To this list, add the Galois extensions:

\medskip
\eject

\noindent {\bf  Lemma 2.5.} {\sl Let $G$ be a finite group.
Any $G$-Galois extension $\psi : R \lr S$ is 
Frobenius. As a Frobenius system one may choose the datum 
$((x_i, y_i)_{i=1, \ldots, m}, \tr )$, where $(x_1, \ldots ,x_m$; $y_1, \ldots , y_m)$
is a Galois basis of $\psi $, and $\tr $ is the trace morphism.
As a Frobenius element one may take the Galois element 
$\displaystyle{ \eta_e  = \sum_{i=1}^m x_i\ot y_i \in S\ot _RS}$.}
\medskip

\noindent {\sl Proof.} This assertion comes as an immediate consequence
from the  definitions and from the following two  equalities, showed in  [27] (D\'emonstration
du lemme 1.4), namely: 
for all $s\in S,$ one has
$$\sum _{i=1}^mx_i\tr (y_is) = \sum _{i=1}^m\tr (sx_i)y_i = s. $$\dm

\medskip
Let $\psi : R \lr S$ be a $G$-Galois extension.
Denote by $C_R(S) = \{ s \in S \ \ \vert \ \ sr = rs, \ \ \forall r \in R\} = S^R$ the centralizer 
of $R$ in $S$.
Clearly the group $G$ acts  on the algebra $C_R(S)$ and  $C_R(S)^G = \ZC (R)$. 
The general theory of Frobenius extensions ([12], 1.3) shows that once 
a Frobenius system $((u_i, v_i)_{i=1, \ldots, m}, \tau)$ is fixed, there is an automorphism 
$\nu : C_R(S) \lr C_R(S)$, called the {\sl Nakayama
automorphism of $\psi $}, which verifies 
$$\tau (\nu (d)s) = \tau (sd),$$
for any $s \in S$ and $d \in C_R(S)$. If the automorphism $\nu$ is inner, the Frobenius extension is 
called {\sl symmetric.}

In the Galois case, taking as Frobenius system
$((x_i, y_i)_{i=1, \ldots, m}, \tr )$, where $(x_1, \ldots ,x_m$; $y_1, \ldots , y_m)$
is a Galois basis of $\psi $, the Nakayama
automorphism is given by the formula
$$\nu (d) = \sum_{i = 1}^m\tr (x_id)y_i = \sum _{g\in G}\sum_{i = 1}^mg(x_id)y_i, $$
for any $d \in C_R(S)$.
If $\psi : R \lr S$ is a centralizing $G$-Galois extension, then $C_R(S) = S,$ and the 
Nakayama automorphism is defined on the whole ring $S.$

\medskip
\noindent {\bf  Lemma 2.6.} {\sl Let $\psi : R \lr S$ be a  $G$-Galois-Azumaya extension such that
the Picard group ${\rm Pic}(R)$ is trivial. Then $S$ is  a symmetric Frobenius $R$-algebra.} 
\medskip

\Dem Denote by ${\rm Inn} (S)$ the group of inner automorphisms
of $S.$ When $R$ has trivial Picard group, the Rosenberg-Zelinsky exact sequence
asserts that ${\rm Inn} (S)  = \Aut (S/R).$ It remains to show that the Nakayama
automorphism leaves $R$ pointwise invariant.
For $r \in R$, one gets 
$$\nu (r) = \sum _{g\in G}\sum_{i = 1}^mg(x_ir)y_i = \sum_{i = 1}^m\bigl(\sum _{g\in G}g(x_i)r\bigr)y_i
= r\sum _{g\in G}\sum_{i = 1}^mg(x_i)y_i = r.$$ \dm

\medskip
\eject
\noindent 2.3. {\sl Galois representations of the braid groups.} 
\smallskip
\noindent
Let $R$ be a commutative algebra and $\psi : R \lr S$ be a Frobenius $R$-algebra.
For any $\displaystyle{ X = \sum _{i = 1}^mx_i \ot y_i \in S\ot _RS}$, define three
elements $X _{12}, X _{23}, X _{13}$ in $S \ot _RS \ot _RS$ by the formulae
$$ X _{12} = \sum _{i = 1}^mx_i \ot y_i \ot 1, \ \ 
X _{23} = \sum _{i = 1}^m1 \ot x_i \ot y_i \ \ {\hbox{\rm and}} \ \
X_{13} = \sum _{i = 1}^mx_i \ot 1 \ot y_i.$$
The properties satisfied by a Frobenius element $\eta$ attached to the extension $\psi : R \lr S$
lead to the fact that $\eta$  satisfies the 
{\sl Frobenius-separability equation} $X_{12}X_{23} = X_{23}X_{13} = X_{13}X_{12}$
[4] (Definition 18), which obviously implies the
Yang-Baxter equation
$X_{12}X_{23}X_{12} = X_{23}X_{12}X_{23}$ (this result appears in [12], Theorem 4.7). 
Therefore, when $\eta $ is invertible in $S\ot _RS$, it
defines a representation of E. Artin's braid group $B_n$ on $n$ strings. Recall that $B_n$ is the group 
generated by symbols $\sigma _1, \ldots , \sigma _{n-1}$
subject to the relations $\sigma _i\sigma _j = \sigma _j \sigma _i$ ($1 \leq i, j \leq n-1$ and $\vert i-j \vert > 1$)
and $\sigma _{i+1}\sigma _i\sigma _{i+1} = \sigma _i\sigma _{i+1}\sigma _i$
($1 \leq i < n-1$).

\medskip
\noindent {\bf  Proposition 2.7.} {\sl For any centralizing
$G$-Galois extension $\psi : R \lr S$, the Galois element 
$\displaystyle{ \eta _e
=}$ $\displaystyle{ \sum _{i = 1}^mx_i \ot y_i}$
verifies the Frobenius-separability equation and therefore the Yang-Baxter equation. 
For any  $G$-Galois-Azumaya extension $\psi : R \lr S$, the Galois element 
$\eta$ is invertible in $S\ot _RS$, hence defines for any integer $n \geq 2$
 a representation $\rho _\psi^n : B_n \lr \Aut (S^{\otimes _{R}n})$ of the braid group
$B_n$.} 
\medskip

\Dem The proposition is  a consequence of Lemma 2.5, of the observation stated above,
and of a crucial result 
by Kadison [12] (Theorem 5.14) asserting that a Frobenius algebra $S/R$ is central separable
if and only if the Frobenius element is invertible in $S\ot _RS$.
Now an invertible  solution $\eta $ of the Yang-Baxter equation 
defines an element $\ell _\eta$ in ${\rm Aut}_R(S \ot S)$, the left multiplication
by $\eta$. This automorphism $\ell_\eta$ is then a solution
of the Yang-Baxter equation for endomorphisms, a so-called 
$\cal R$-matrix, 
hence gives rise, for any
$n \geq 2$, to a representation $\rho _\psi^n$ of the braid group $B_n$ 
in $S^{\otimes  _{R}n}$ given by $\rho _\psi^n (\sigma _i ) = \id _S^{\otimes (i-1)}
\ot \ell_\eta \ot \id _S^{\otimes (n-i-1)}$ [14] (Corollary X.6.9). 
\dm

\medskip
\noindent {\bf  Definition 2.8.} For any $G$-Galois-Azumaya extension $\psi : R \lr S$,  the representation 
$\rho _\psi^n$ is called the {\sl Galois representation of $B_n$ attached to $\psi$.}

\medskip 
\noindent {\bf  Corollary 2.9.} 
{\sl 
A morphism $(\varphi , \tilde \varphi )$  from $(\psi : R \lr S)$ to $(\psi ' : R' \lr S')$
in ${\goth GalAzum}(G)$ induces, for any integer $n\geq 2$, a morphism $\rho _\psi^n \lr \rho _{\psi'}^n$
of Galois representations of the braid group $B_n$.}

\medskip
\noindent {\sl Proof.} Fix $n\geq 2$.  
From Corollary 1.6, one deduces that the diagram
$$\diagram {
S^{\otimes _Rn} &\hfll{\rho _\psi^n}{} &   S^{\otimes_R n} \cr
\vfl {\varphi ^{\otimes n}}{} &  & \vfl {}{\varphi ^{\otimes n}}&\cr
S'^{\otimes _{R'}n} &\hfll{}{\rho _{\psi'}^n}& 
S'^{\otimes _{R'}n} &\cr } $$
is commutative.\dm

\medskip
\noindent {\sl  Example 2.10.} 
Let ${\displaystyle H = ({{a,b}\over F})}$ be a quaternion algebra.
In [27] (\pa \ 2.3), we have seen 
that $$\eta = {{1}\over 4}(1\ot 1 + {{i\ot i}\over a} + {{j\ot j}\over b}
- {{k\ot k}\over ab})$$
is the Galois element of the $V$-Galois-Azumaya extension
$\psi : F \lr H.$ Let us mention on the way that the
Nakayama automorphism $\nu $ is in this case equal to the identity of $H$. Indeed, the trace map
$\tr : H \lr F$ is given by $\tr (c_0 + c_1i + c_2j + c_3k) = 4c_0,$ where $c_0 + c_1i + c_2j + c_3k \in H$.
Take as Galois basis $x_1 = 1$, $x_2 = i$, $x_3 = j$, $x_4 = k$~; $y_1 = 1/4$,  $y_2 = i/4a$,
$y_3 = j/4b$, $y_4 = -k/4ab$. Then $\displaystyle{\nu (c) = \sum_{i = 1}^m\tr (x_ic)y_i = c}$.

Let us now describe the automorphism 
$\ell_\eta$ with more details. Decompose the 16-dimensional $F$-vector space $H \ot _FH$ into four
subspaces $V_1, V_i, V_j, V_k$ each of them of dimension 4 in the following way.
$$\vbox{\catcode`\*=\active \def*{\hphantom{0}}
\offinterlineskip
\ialign{
\vrule height11pt depth4pt \vrule\quad#\hfil\unskip\hfil
\quad&\vrule\quad#\hfil\unskip\hfil\quad\vrule
&&\quad#\unskip\hfil\quad\vrule\cr
\noalign{\hrule}
\ \ Vector space\ \ & Basis \cr
\noalign{\hrule}
\ \ \ $V_1$ \ \ &   $(1 \ot 1, i \ot i, j \ot j, k \ot k)$  \cr
\ \ \ $V_i$ \ \ &  $(1 \ot i, i \ot 1, j \ot k, k \ot j)$ \cr
\ \ \ $V_j$ \ \ & $(1 \ot j, j \ot 1, k \ot i, i \ot k)$ \cr
\ \ \ $V_k$ \ \ &  $(1 \ot k, k \ot 1, i \ot j, j \ot i)$\cr
\noalign{\hrule}}}$$
Then $H \ot _FH =  V_1 \oplus V_i \oplus V_j\oplus V_k$ and  each of the subspaces $V_1,
V_i, V_j, V_k$ is stable under $\ell_\eta$.
In the given bases, $\ell_\eta$ is represented by the matrix $\Xi _u$ in $V_u$
($u = 1, i, j, k$), where
$$\matrix{\Xi _1 & = & {1\over 4}\pmatrix 
{1 & a & b & -ab \cr
{1\over a} & 1 & {-1\over ab} & {1\over b} \cr
{1\over b} & {-1\over ab} & 1 & {1\over a} \cr
{-1\over ab} & {1\over b} & {1\over a} & 1 \cr}, &{\rm } 
&\Xi _i &= & {1\over 4}\pmatrix 
{1 & 1 & 1 & {-1\over ab}  \cr
1 & 1 & {-1\over ab} & 1 \cr
{-1\over b} & {-1\over ab} & 1 & {-1\over a} \cr
{-1\over ab} & {-1\over b} & {-1\over a} & 1 \cr}, \cr
&{} &{} & {} & {} & & \cr
&{} &{} & {} & {} & & \cr
\Xi _j & = & {1\over 4}\pmatrix 
{1 & 1 & 1 & -1 \cr
1 & 1 & -1 & {1\over ab} \cr
{1\over ab} & {1\over a} & 1 & {-1\over b} \cr
{1\over a} & {1\over ab} & {-1\over b} & 1 \cr}, &{\rm and} 
&\Xi _k & = & {1\over 4}\pmatrix 
{1 & 1 & 1 & -1  \cr
1 & 1 & -1 & 1 \cr
{-1\over a} & {1\over b} & 1 & {1\over ab} \cr
{1\over b} & {-1\over a} & {1\over ab} & 1 \cr}.\cr}$$

The quaternionic extension ${\displaystyle F \lr  H}$
gives moreover rise to a representation of the symmetric groups ${\goth S}_n$ ($n\geq 2$).
Indeed the square of $\eta $ in the $F$-algebra
$H \ot _FH$ is
$$\eta ^2 = {1 \over 4}(1 \ot 1).$$
Thus $(2\eta )^2 = 1\ot 1$, that is to say that $\ell _{2\eta}$ is involutive. 
The assertion follows from Moore's presentation of the symmetric groups:
for any $n \geq 2$, ${\goth S}_n$ is the group 
generated by symbols $\sigma _1, \ldots , \sigma _{n -1}$
subject to the relations
$\sigma _i\sigma _j = \sigma _j \sigma _i$ ($1 \leq i, j \leq n-1$ and $\vert i-j \vert > 1$),
$\sigma _{i+1}\sigma _i\sigma _{i+1} = \sigma _i\sigma _{i+1}\sigma _i$
($1 \leq i < n-1$) and  $\sigma _i ^2 = 1$ 
($1 \leq i \leq n-1$).

\medskip

\noindent 2.4. {\sl A fundamental example of a Galois-Azumaya extension: 
the $n$'th power norm residue algebra.} 
\smallskip
\noindent
Let $F$ be a commutative field and $n \geq 2$ an integer. 
Suppose that the multiplicative group $F^{\lcross}$ of $F$ contains a primitive $n$-th root of
unity $\zeta$, that is an element of order $n$ in $F^{\lcross}$. Choose once and for all two elements
$a, b \in F^{\lcross}$. With these datas, following
Milnor [24] (\pa \ 15) (we refer  also to [8] (\pa \  11)
or to [29] (15.4)), one explicitly constructs  a central simple associative algebra in which $n$ plays the same r\^ole
as $2$ does for quaternion algebras. This is done in the following way (see [21]  (16.22, iv)):

Let $Q$ be the quotient ring $F[X]/(X^n - a)$ of the polynomial ring $F[X]$
in one indeterminate $X$ by the ideal $(X^n -a)$. Denote by $x$ the image of
$X$ by the canonical projection $F[X] \lr Q$. By the euclidian division algorithm, 
$\{ 1, x, x^2, \ldots , x^{n -1}\}$ is a basis of the $F$-vector space $Q$. 
The $F$-algebra homomorphism $s : F[X] \lr F[X]$ defined by $s(X) = \zeta X$
induces a $F$-algebra automorphism $\sigma : Q \lr Q$ that is of order $n$
in ${\rm Aut}(Q)$. Identify the cyclic group $\langle \sigma \rangle$ generated by 
$\sigma$ with $\ZB /n\ZB$.
Take now the free $Q$-module $S$ based on the set $\{ u^m \ \ \vert \ \ m = 0, \ldots, n-1\}$
indexed by the group $\ZB /n\ZB$.
As an $F$-vector space, $S$ has  the basis $\{ x^mu^{m'} \ \ \vert \ \ m, m' = 0, \ldots, n-1\}$.
Put on $S$ the multiplication determined by the relations
$$x ^n = a, \ \ \, u^n = b, \ \ \ {\rm and} \ \ \ ux = \zeta xu.$$
The $F$-vector space $S$ becomes
an algebra $S$, called the {\sl $n$'th power norm residue algebra} (or a {\sl symbol algebra} [16])
and denoted by $(a, b, \zeta )_F$. It is an Azumaya algebra of degree $n$ ([8], \pa \  11).

\medskip
\noindent {\sl Remark:} The $n$'th power norm residue algebras clearly generalize quaternion 
algebras: let $F$ be a field of  characteristic different from 2 and  $a, b \in F^{\lcross}$. Then
 $\displaystyle{({{a,b}\over F})}$ is isomorphic to $(a, b, -1)_{F}$ [29].
Therefore, following Pierce, one denotes $(a, b, \zeta )_F$
also by $\displaystyle{({{a,b}\over {F, \zeta}})}$.

\medskip

\noindent {\bf  Theorem 2.11.} {\sl Let $F$ be a commutative  field.
Take two elements $a$ and $b$
in $F^{\lcross}$, and  $n \geq 2$ an integer. 
Suppose that the multiplicative group $F^{\lcross}$ of $F$ has an
element $\zeta$ of order $n$.  Let $(a, b, \zeta )_F$
be the $n$'th power norm residue algebra. Then the extension $F \lr (a, b, \zeta )_F$
is  $(\ZB /n\ZB)^2$-Galois-Azumaya with Galois element
$$\eta = \sum _{r, s = 0}^{n-1} {{\zeta ^{rs}} \over {abn^2}} x^ru^s \ot x^{n-r}u^{n-s}.$$}

\Dem Set $S = (a, b, \zeta )_F$ and $G = (\ZB /n\ZB)^2$.
It is already  known that $S$ is an Azumaya $F$-algebra
([21], Proposition 16.24). It remains to show that
the extension $F \lr S$
is $G$-Galois. 

Denote by $(i, j)$, with $0 \leq i,j \leq n-1$, the elements of 
$G$. It is easy to  see that the formula
$$(i,j)\cdot (qu^m) = \sigma ^j(q)\zeta ^{im}u^m$$
defines an action of $G$ on 
$S$ (here $0 \leq i,j, m \leq n-1$ and $ q \in Q$). 
Thus, if $\alpha = (1,0)$ and $\beta = (0,1)$ are the two 
canonical generators
of G, one gets $\alpha (x) = x$, $\alpha (u) = \zeta u$, $\beta (x) = \zeta x$, and $\beta (u) = u.$
Suppose now that ${\displaystyle \sum _{m = 0}^{n-1}q_mu^m \in S^G}.$
Then, in particular ${\displaystyle  (0, 1)\cdot \sum _{m = 0}^{n-1}q_mu^m =
\sum _{m = 0}^{n-1}q_mu^m}$, that is
${\displaystyle  \sum _{m = 0}^{n-1}\sigma (q_m)u^m =
\sum _{m = 0}^{n-1}q_mu^m}$. But $(u^m)_{0 \leq m \leq n-1}$ is a $Q$-basis
of $S$,
thus, for any $m$, one has $\sigma (q_m) = q_m$.

Pick now an element ${\displaystyle q = \sum _{k = 0}^{n-1}a_kx^k \in Q}$
which is invariant under   $\sigma$.
This means that ${\displaystyle \sum _{k = 0}^{n-1}a_kx^k = \sum _{k = 0}^{n-1}a_k\zeta^kx^k}$.
Again $(x^k)_{0 \leq k \leq n-1}$ is an $F$-basis of $Q$, therefore 
$a_k = a_k\zeta^k$, for all $k$, thus $a_k = 0$ as soon as $k \geq 1$.
So $q$ belongs to $F.$ Conversely, it is clear that an element $a_0 \in F$
 is invariant under $G.$ Hence $S^G = F.$

Before we end the proof, we write down an important formula that is readily
deduced from the multiplication law in $S$ given above.
For any $r,s,i,j,k,l \in \{0, \ldots , n-1\}$, an easy computation gives 
$$x^ru^s(i,j)(x^ku^l) = \zeta ^{il + jk  + ks}x^{r+k}u^{s+l}.$$
In particular, 
$$x^ru^s(i,j)(x^{n-r}u^{n-s}) = \zeta ^{-(is + jr + rs)}ab. \eqno (2)$$
Having this rule of calculus in view, we claim that the following equality is true:
$$\Gamma \bigl(\sum _{r, s = 0}^{n-1} {{\zeta ^{rs}} \over {abn^2}} x^ru^s \ot x^{n-r}u^{n-s}\bigr) = 
\delta _{(0,0)}. \eqno (3)$$ This is then sufficient to show that $\psi : F \lr S$
is $G$-Galois, since formula $(3)$ exactly means that there exists a Galois basis (see
[27], D\'efinition 1.6 and Th\'eor\`eme 1.5). Let us  prove $(3)$.

 According to $(2)$, the left-hand side is equal to
${\displaystyle {{1} \over {n^2}}
\sum _{(i,j) \in G}\sum _{r,s = 0}^{n-1}\zeta ^{-(jr + is)}\delta _{(i,j)}.}$ 
The coefficient of 
$\delta _{(i,j)}$ is
$\displaystyle{\sum _{r,s = 0}^{n-1}\zeta ^{-(jr + is)} =  
\sum _{r = 0}^{n-1}\zeta^{-jr}(\sum _{s = 0}^{n-1}(\zeta ^{-i})^s)}$. 
But the cyclotomic sum $\displaystyle{\sum _{t = 0}^{n-1}(\zeta ^{-k})^t}$ is equal either to $0$, when
$k$ is different from $0$, or to $n$, when  $k = 0$.
\dm

\medskip
\noindent {\sl Remark:} Taking $n = 2$, $x = i$, $u = j$, $xu = k$, one verifies that
the tensor decomposition of the Galois element in  Theorem 2.11 coincides with 
 that of the $V$-Galois-Azumaya extension
$\psi : F \lr H$ given in Example 2.10.
\medskip
Similarly to the case of Hamilton's quaternion, where Theorem 1.3
shows that the inclusions $\RB \lr \CB$ and $\CB \lr \HB$ are both $\ZB/2\ZB$-Galois extensions,
 one deduces, with the notations adopted above, the following result:

\medskip
\noindent {\bf  Corollary 2.12.} {\sl Let $F$ be a commutative  field.
Take two elements $a$ and $b$ 
in $F^{\lcross}$, and  $n \geq 2$ an integer.
Suppose that the multiplicative group $F^{\lcross}$ of $F$ has an
element $\zeta$ of order $n$.  
Let $Q$ be the quotient ring $F[X]/(X^n - a)$. Then the extensions $F \lr Q$
and $Q \lr (a, b, \zeta )_F$
are $\ZB /n\ZB$-Galois. Moreover $Q$ is the maximal commutative subring of $(a, b, \zeta )_F$.} 
\medskip
The maximality of $Q$ in $(a, b, \zeta )_F$ is shown in  [21] (Proposition 16.24).
The definition of $Q$ is in fact the usual cyclic extension performed over an $n$-kummerian
ring ([10], Chapter 0, \pa \ 5) (recall that an {\sl $n$-kummerian
ring} is a commutative ring which contains $n^{-1}$
and a root of the $n$-cyclotomic polynomial). 

\medskip
\noindent {\bf  Corollary 2.13.} {\sl Let $n \geq 2$ be an integer
and $F$ be a commutative field such the multiplicative group $F^{\lcross}$ of $F$ has an
element $\zeta$ of order $n$.  
Then the extension $F \lr M_n(F)$ is 
$(\ZB /n\ZB)^2$-Galois-Azumaya.} 
\medskip
\Dem The matrix algebra $M_n(F)$ is a particular case of $n$'th power norm residue algebras:
choose a primitive $n$-th root of
unity $\zeta$ and $a, b \in F^{\lcross}$, then
$M_n(F)$ is isomorphic to $(a, b, \zeta )_F$ if (and only if) 
$b \in N_{E/F}(E^{\lcross})$ or equivalently 
if (and only if) 
$a \in N_{K/F}(K^{\lcross})$, where $E = F(a^{1 \over n})$, $K = F(b^{1 \over n})$, and
$N_{L/F}$ is the norm map $L^{\lcross} \lr F^{\lcross}$ for any finite extension $L/F$ ([21], Chapter 16).
For instance  $(1, b, \zeta )_F \cong M_n(F).$\dm

\bigskip
\noindent {\bf 3. The Brauer-Galois group of a commutative ring.}
\smallskip
\noindent 3.1. {\sl  The category ${\goth Gal}_R$.} 
\smallskip
\noindent
We study some
properties of the  two subcategories of ${\goth Gal}$ 
obtained by fixing a Galois group, or by fixing simultaneously  a base ring
and a Galois group (of course, an analogue work can be done starting
with one of the categories ${\goth Galstr}$, ${\goth Galcent}$ or ${\goth Galcom}$).

\medskip
The functor 
$\pi _{alg} : {\goth Gal} \lr {\goth Alg}_k$,
given on objects by $\pi _{alg}(G, \psi : R \lr S) = R$, enables us to consider the fibre 
category over a fixed algebra  $R$. Denote it by
${\goth Gal}_R$. Its objects are the couples $(G, \psi : R \lr S)$; a
morphism from $(G, \psi : R \lr S)$ to $(G', \psi ' : R \lr S')$
in ${\goth Gal}_R$
is a couple $(f, \varphi)$, where $f : G'\lr G$ is a morphism of groups
and $\varphi : S \lr S'$ a morphism of algebras. These data are subject to the two conditions: 
$\varphi  \circ \psi  = \psi ' $ and
$\varphi  \circ f(g') = g' \circ \varphi ,$ for any $g'\in G'$.
Necessarily, $\varphi $ is a morphism of $R$-rings.

\medskip

\noindent {\bf  Theorem 3.1.} {\sl Let
$R$ be a commutative algebra. Suppose that $(G, \psi : R \lr S)$ and 
$(G', \psi ' : R \lr S')$ are two objects in ${\goth Galcent}_R$
(respectively in ${\goth GalAzum}_R$). Then
the canonical extension $\psi'' : R \lr S\ot _RS'$ is a 
centralizing (respectively central) $G \cross G'$-Galois extension that we
 denote  by $\psi \ot \psi '$.}
\medskip

\noindent {\sl Proof.} Let $(G, \psi : R \lr S)$ and 
$(G', \psi ' : R \lr S')$ be two centralizing Galois extensions.
Since the ring $R$ is commutative,
the tensor product $S\ot _RS'$ comes with the multiplication given by 
$(s\ot s').(t\ot t') = st\ot s't'$, with $s,t \in S$ and $s',t' \in S'$.
The group $K = G\cross G'$ acts on $U = S\ot _RS'$ by 
$(g, g')\cdot (s\ot s')
= g(s)\ot g'(s'),$ for $(g , g') \in K$ and $s \ot s' \in U.$
We must show that $\psi '' : R \lr U$ is  $K$-Galois. In order to do that, we successively
prove that $R$ can be identified with the ring $U^K$ of invariants of $U$ under $K$, 
that the associated morphism $\Gamma : U\ot _RU \lr U(K)$ is an  isomorphism, 
and finally that the extension $R\lr U$ is faithfully flat.

1) The inclusion $R \subseteq U^K$ is obvious. In order to prove 
$U^K \subseteq R$, we use an argument inspired by
Chase-Harrison-Rosenberg [5] (Lemma 1.7) which involves descent.
By [27] (Lemme 1.4.2), the trace morphisms
$\tr _{S/R}$ and $\tr _{S'/R}$ are surjective. Choose $s_0\in S$
and $s'_0\in S'$ with $\tr _{S/R}(s_0) = \tr _{S'/R}(s'_0) = 1$. For any $u \in U^K$, the following equalities hold

$$\eqalign{u &= \bigl((\tr _{S/R} \ot \tr _{S'/R} )(s_0 \ot s'_0)\bigr).u 
=  \bigl(\sum_{g \in G}g(s_0) \ot \sum_{g' \in G'}g' (s'_0)\bigr).u \cr
& = \sum_{(g,g') \in K}\bigl(g(s_0) \ot g'(s'_0)\bigr).u  = \sum_{(g,g') \in K}(g , g' )
\bigl((s_0 \ot s'_0).u\bigr) = (\tr _{S/R} \ot \tr _{S'/R} )\bigl((s_0 \ot s'_0).u\bigr),\cr}$$
that belongs to ${\rm Im}(\tr _{S/R} \ot \tr _{S'/R} )$ = $R\ot _RR \cong R. $

 2) The morphism $\Gamma : U\ot _RU \lr U(K)$ defined by
$$\Gamma \bigl((s\ot s')\ot (t\ot t')\bigr) =
\sum_{(g,g') \in K}\bigl(sg(s')\ot tg'(t')\bigr)\delta_{(g,g')}$$
is an isomorphism of left $R$-modules. Indeed, one may easily see that
$\Gamma$ is the composition \goodbreak \noindent $\beta \circ (\Gamma _S \ot \Gamma _{S'}) \circ \alpha$,
where $\alpha : U\ot _RU \lr (S\ot _RS)\ot _R(S'\ot _RS')$ 
and $\beta : S(G)\ot _RS'(G') \lr$ $ (S\ot _RS')(K)$ are the canonic isomorphisms.

 3) Proposition 1.2 in [27] asserts that when  $U^K = R$ and when
$\Gamma : U\ot _RU \lr U(K)$ is an isomorphism, then the
extension $R\lr U$ is faithfully flat if and only if
$R$ is a direct summand of the (right or left) $R$-module  $U$. 
The extension $R\lr S$ (respectively $R\lr S'$) being Galois, $R$ is a direct summand of the $R$-module 
$S$ (respectively $S'$). Hence $S\ot _RS'$ contains $R\ot _RR \cong R$ as a direct summand, because
the tensor product preserves direct limits and consequently also direct sums.

Suppose now that $(G, \psi : R \lr S)$ and 
$(G', \psi ' : R \lr S')$ are Galois-Azumaya. Since the centre
$\ZC (S \ot _R S')$ of $S\ot _RS'$ is isomorphic to $\ZC (S) \ot _R \ZC (S')$
[7], the extension $\psi'' : R \lr S\ot _RS'$ is  $G \cross G'$-Galois-Azumaya. \dm

\medskip

\noindent {\sl Examples.} 

\noindent 1.-- Let $F$ be a commutative field
of characteristic different from 2 and  let
$a _m, b_m$ ($m = 1, \ldots , n$) be $2n$ non zero elements in $F$.
Denote by $H_m$ the quaternion algebra
${\displaystyle ({{a_m,b_m}\over F})}$.
By Theorem 3.1, the $n$-fold tensor product extension $F \lr H_1\otimes _F \ldots \otimes _F H_n$
is $V^n$-Galois. Observe that the tensor product
$H_1\otimes _F  H_2$ of two quaternion algebras may be not a skew field:
this is the case if (and only if) $H_1$ and $H_2$ 
have a common splitting field which is separable quadratic
over $F$ [8] (\pa \  14, Theorem 6).
For instance, let $F$ be the field  $\RB$ of the real numbers
and  $a = b = -1$. Thus ${\displaystyle ({{-1,-1}\over \RB})}$ 
is the field  
$\HB$ of Hamilton's quaternions and the algebra $\displaystyle{\HB\ot_{\RBIc}\HB }$ is not a skew field.
\smallskip

\noindent 2.-- Let $F$ be a commutative field
of characteristic different from 2 and $A$  central simple algebra of degree 4 and exponent 2
(recall that the exponent of an $F$-Azumaya	algebra is the order, in the sense of Group Theory,
of its class $[A]$ in the Brauer group ${\rm Br}(F)$).
Then the extension $F \lr A$ is $V^2$-Galois-Azumaya. Indeed, following a
well-known result of Albert [16], $A$ is necessarily
a biquaternion algebra, that is a tensor product of two quaternion algebras over $F$.

\medskip

\noindent {\bf  Proposition 3.2.} {\sl Fix
a commutative algebra $R$. Let $\psi : R \lr S$ (respectively $\psi ' : R \lr S'$)
be a centralizing Galois extension with Galois group $G$ (respectively $G'$)
of order $n$ (respectively $n'$). Assume that the integer $nn'$ is invertible in $S\ot _RS'$.
Then the extension $\varepsilon _1 : S \lr S \ot _RS'$ given by $\varepsilon _1 (s) = s\ot 1$ $(s\in S)$
is $G$-Galois and the extension $\varepsilon _2 : S' \lr S \ot _RS'$ given by $\varepsilon _2 (s') = 1\ot s'$ $(s'\in S')$
is $G'$-Galois.}

\medskip

\noindent {\sl Proof.} Because of the condition on $nn'$, the
$G \cross G'$-Galois extension $\psi \ot \psi ' : R \lr S\ot _RS'$ becomes strict. Take 
 the subgroup
$H = \{1\} \cross G'$ of $G \cross G'$. By Theorem 1.3, 
the inclusion $\theta : U = (S\ot _RS')^{H} \lr S\ot _RS'$ is $H$-Galois. It is enough to show
that $(S\ot _RS')^{H} \cong S.$

Consider the right $S'$-module $M' = S \ot _RS'$. Identifying $\{1\} \cross G'$ with $G'$, the group action
on $M'$ is given by $g'(s\ot s') = s \ot g'(s')$ ($g'\in G', s\in S, s'\in S'$). So
$M'$ is  a $G'$-Galois module. By Galois descent for modules [26] (Corollary 4.16), 
the module $N' = {{M'}}^{G'}$ is such that $N' \ot_RS' \cong M'$, that is
$(S \ot _RS)^{G'} \ot_RS' \cong S\ot _RS'$. By faithfully flatness of $\psi '$,
one concludes that ${(S \ot _RS)}^{G'} \cong S$ and that  $\varepsilon _2 : S' \lr S \ot _RS'$ 
is $G'$-Galois.
Similarly $\varepsilon _1 : S \lr S \ot _RS'$ 
is $G$-Galois. \dm

\medskip

\noindent {\sl Example.} Take 
 $F \lr H$ 
a quaternionic extension over a field $F$ of characteristic different from 2.  By
Proposition 3.2, the two extensions $\varepsilon _i : H \lr H \ot _FH$ $(i=1,2)$ are
both $V$-Galois.

\medskip

\noindent 3.2. {\sl  The Brauer-Galois group.}
\smallskip
\noindent {\bf  Corollary  3.3.} {\sl Let $R$ be a commutative algebra. Denote by
${\rm Br}_{\rm Gal}(R)$ the subset of the Brauer group of $R$ 
consisting of those
classes $\xi \in {\rm Br}(R)$ for which there exists a finite group $G$ and 
a $G$-Galois-Azumaya extension $\psi : R \lr S$, such that $\xi$ is equal to the
class $[S]$ of $S$ modulo stable isomorphisms. Then ${\rm Br}_{\rm Gal}(R)$
is a subgroup of ${\rm Br}(R)$.}

\medskip 

This abelian group ${\rm Br}_{\rm Gal}(R)$ is called the {\sl Brauer-Galois group of $R$}.

\medskip
\noindent {\sl Proof.} The neutral element of ${\rm Br}(R)$ is the 
class $[R]$. Since ${\rm id} : R \lr R$ is a
$\{1\}$-Galois-Azumaya extension, $[R]$ belongs 
to  ${\rm Br}_{\rm Gal}(R)$, and this set is therefore non-empty.

Let $\xi_1, \xi_2\in {\rm Br}(R)$. Represent these classes
by $\psi _1 : R \lr S_1$ and $\psi _2 : R \lr S_2$, two
Galois-Azumaya extensions with Galois groups $G_1$ and $G_2$. The product 
 $\xi_1 . \xi_2\in {\rm Br}(R)$ can then be represented by the $R$-Azumaya algebra $S_1 \ot _RS_2$
(the multiplication in $S_1 \ot _RS_2$ being given by 
$(s_1 \ot s_2)\cdot (t_1 \ot t_2) = s_1t_1\ot s_2 t_2$,
for any $s_1,t_1 \in S_1$ and $s_2,t_2 \in S_2$, is well defined since $R$ is commutative
and $\psi _i$ ($i = 1,2$) is central).
By Theorem 3.1, $\psi _1 \ot  \psi _2 : R \lr S_1 \ot _R S_2$ is a $G_1 \cross G_2$-Galois
extension. Hence the set ${\rm Br}_{\rm Gal}(R)$ is stable under the product in ${\rm Br}(R)$.

The symmetric element of $[S]$ can be represented
by the opposite algebra $S^{\rm o}$. But if $\psi : R \lr S$ is a $G$-Galois-Azumaya extension,
the group $G$ acts on $S^{\rm o}$ by $g(s^{\rm o}) = \bigl(g(s)\bigr)^{\rm o}$. Since $R$ is
the centre of $S$, the map $\psi ^{\rm o} : R \lr S^{\rm o}$ is well defined. It is a
$G$-Galois extension by Lemma 1.4.
\dm
\medskip

\noindent {\sl Remark:} A similar construction can be provided in the Hopf-Galois context when one replaces
Galois-Azumaya extensions by Hopf-Galois-Azumaya extensions. This gives rise to
a group ${\rm Br}_{\rm HopfGal}(R)$, the {\sl Brauer-Hopf-Galois group of $R$}, which should be
less interesting than ${\rm Br}_{\rm Gal}(R)$ since it is bigger: the inclusions
${\rm Br}_{\rm Gal}(R) \subseteq {\rm Br}_{\rm HopfGal}(R) \subseteq {\rm Br}(R)$ hold.

\medskip
\noindent 3.3. {\sl  Some straightforward facts about the Brauer-Galois group.}  

\smallskip
\noindent 1) The Brauer-Galois group ${\rm Br}_{\rm Gal}(S)$ is obviously trivial when
${\rm Br}(S)$ is trivial. Therefore the group
${\rm Br}_{\rm Gal}(S)$ is trivial for any finite field $S = \FB _q$ (by Wedderburn's Theorem
on finite division rings), as well as for any  algebraic extension of 
 $\FB_q$. It is also trivial for any algebraically closed  field
 or for any field  of transcendence degree one  over an algebraically closed  field (Tsen's Theorem) [9]. 
Other examples of commutative fields with trivial Brauer group can be found in [31] (p. 170). 
\medskip
\noindent 2) Let $F$ be a commutative field
of characteristic different from 2 and  
$a, b \in F^{\lcross}$. Since ${\displaystyle ({{a,b}\over F})^{\rm o} \cong ({{a,b}\over F})}$, 
the quaternion algebra ${\displaystyle ({{a,b}\over F})}$ is 
an element of order at most $2$ in ${\rm Br}_{\rm Gal}(F)$. It is exactly of order two 
if and only if
the Hilbert symbol $(a,b)_F$ is equal to $1$, in other words 
if and only if ${\displaystyle ({{a,b}\over F})}$ is not isomorphic 
to the matrix algebra 
$M_2(F)$ [30].
\medskip
\noindent 3)   For any commutative field $F$ 
of  characteristic different from 2, the Brauer-Galois group ${\rm Br}_{\rm Gal}(F)$ 
contains the quaternion group
${\rm Quat}(F)$, that is the subgroup of ${\rm Br}(F)$ generated by the quaternion algebras
(these lie in the $2$-torsion).
Under suitable conditions on the characteristic of $F$ (see Theorem 2.11),
the Brauer-Galois group ${\rm Br}_{\rm Gal}(F)$ 
contains the subgroup of ${\rm Br}(F)$ generated by the $n$'th power norm residue algebras
(these lie in the $n$-torsion).

\medskip
Next we shall use a deep 
result due to Merkurjev and Suslin
in order to discuss the relationship between the $n$-torsion
of ${\rm Br}_{\rm Gal}(F)$  and 
the $n$'th power norm residue algebras, for any $n \geq 2$.
Before doing that, we have to make a short digression  into $K$-theory.

\medskip
\noindent 3.4. {\sl  Connection with Milnor's $K_2$ and the main Theorem.}
\smallskip
\noindent
For any abelian group $A$  and any integer $n \geq 2$, 
denote by $_nA$ the $n$-torsion subgroup of $A$, that is the part 
of $A$ annihilated by the multiplication
by $n$. 
Let $F$ be a commutative field. Following a famous result of Matsumoto
[21] (16.49), the Milnor's $K_2$-group $K_2(F)$
is generated by the symbols $\{ a, b \}$ (with $a, b \in F^{\lcross}$) subject to the relations
$$\{ ab, c \} = \{ a, c \} \ \{ b, c \}, \ \ \{ a, bc \} = \{ a, b \} \ \{ a, c \}, \ \ {\hbox {\rm and}} \ \ 
\{ a, b \} = 1 \ \ {\hbox {\rm if}} \ \ a + b = 1.$$

\noindent Suppose that the field $F$ is such that the multiplicative group $F^{\lcross}$
has an element $\zeta$ of finite order $n$. Then 
the assignment  $F^{\lcross}\ot _{\ZBIc} F^{\lcross} 
\lr {\rm Br}(F)$ defined by $r_{n,F} (\{ a , b \}) =
[(a, b, \zeta )_F]$ is a Steinberg symbol, hence defines a homomorphism
of groups  $R_{n,F} : K_2(F) 
\lr {\rm Br}(F)$. Furthermore $R_{n,F}$ annihilates $_nK_2(F)$ and takes its values in 
$_n{\rm Br}(F)$ (for all these facts, see [24]). Thus there exists a homomorphism 
$${\bar R_{n, F}} : K_2(F)/nK_2(F) = K_2(F)\ot _{\ZBIc}(\ZB/n\ZB) \lr _n{\rm Br}(F),$$ 
called the {\sl Galois symbol},
which is, following a celebrated and difficult Theorem 
proved by Merkurjev and Suslin [23], an isomorphism of groups (for the formulation we use, see [21]).

\medskip
\noindent {\bf  Theorem 3.4.}
{\sl Let $F$ be a commutative  field, 
and $n \geq 2$ an integer. 
Suppose that the multiplicative group $F^{\lcross}$ of $F$ has an
element $\zeta$ of order $n$. 
Then the  Steinberg symbol $F^{\lcross}\ot _{\ZBIc} F^{\lcross} 
\lr {\rm Br}(F)$ factorizes through ${\rm Br}_{\rm Gal}(F)$, hence 
$$_n{\rm Br}_{\rm Gal}(F) = {_n{\rm Br}}(F).$$

\noindent In particular, if $F$ is a commutative field
of  characteristic different from 2, then $_2{\rm Br}_{\rm Gal}(F) = {_2{\rm Br}}(F) = {\rm Quat}(F).$}

\medskip

\Dem This result comes readily from Theorem 2.11, Corollary 3.3 and the Theorem of Merkurjev-Suslin. 
\dm
\medskip

\noindent The case $n = 2$
already follows from a theorem of Merkurjev [22], which shows that 
the $2$-torsion of the Brauer group is ${\rm Quat}(F)$.
This result implies in particular that

\item {--} The Brauer-Galois group ${\rm Br}_{\rm Gal}(\RB )$ of the real numbers is equal to $\ZB/2\ZB$,
since ${\rm Br}(\RB )$ is equal to $\ZB/2\ZB$ and is generated by Hamilton's quaternions $\HB$ [28] (Example 3.8);

\item {--} The Brauer-Galois group ${\rm Br}_{\rm Gal}(\QB )$ of the rational numbers is infinite,
since ${\rm Quat}(\QB)$ is infinite  (see the classification of quaternion algebras via quadratic forms [9]).

\medskip
\noindent {\bf  Corollary 3.5.} {\sl Let $F$ be a commutative  field
of characteristic zero such that, for any $n \geq 2$, 
the multiplicative group $F^{\lcross}$ of $F$ has an
element of order $n$. Then 
$${\rm Br}_{\rm Gal}(F) = {\rm Br}(F).$$} 

\Dem The group ${\rm Br}(F)$ is torsion ([9], Corollary 4.15), hence as sets
$${\displaystyle {\rm Br}(F) = \bigcup _{n \geq 1}{}{_n{\rm Br}}(F) = 
\bigcup _{n \geq 1}{}{_n{\rm Br}}_{\rm Gal}(F)
= {\rm Br}_{\rm Gal}(F).}$$ \dm

\medskip
\noindent {\sl  Remarks.}

\noindent 1) Observe that the arithmetic conditions in Corollary 3.5 are sufficient
but not necessary since ${\rm Br}_{\rm Gal}(\RB)$ = ${\rm Br}(\RB)$ = $\{ [\RB ], [\HB ]\}$.

\smallskip
\noindent 2) Notice that the conditions that have to be satisfied by the field
$F$ do not imply that $F$ is algebraically closed (in which case
Corollary 3.5 would be vacuous, the Brauer group of an algebraically closed field being trivial). 
For example take $\displaystyle{F = \QB _{\rm ab} = \bigcup _{n \geq 1} 
\QB (\mmu _n (\CB ))}$, the field
obtained by adjoining all roots of unity to the rationals. Then 
$F$ is of characteristic zero and is strictly contained in the algebraic closure 
$\bar {\QB}$ of $\QB$. Indeed, in the 
absolute Galois group ${\rm Gal}(\bar {\QB}/\QB)$,  the only non-trivial torsion elements 
have order $2$ ([19], VI, \pa \  9, Example 1), whereas  $\displaystyle{{\rm Gal}(\QB _{\rm ab}/\QB)
 = \limproj _{n}{\rm Gal}\bigl(\QB ( \mmu _n (\CB ))/\QB\bigr)
= \limproj _{n}(\ZB/n\ZB)^{\lcross}}$ is 
isomorphic to the group $\hat {\ZB}^{\lcross}$ of invertible elements 
 of the Pr\"ufer ring $\displaystyle{\hat {\ZB} = 
\limproj _{n}(\ZB/n\ZB) = \prod _{p \  {\hbox {\sevenrom prime}}}\ZB_p}$; therefore
$\displaystyle{{\rm Gal}(\QB _{\rm ab}/\QB) \cong \bigl(\prod _{p \not= 2\ {\hbox {\sevenrom prime}}}\ZB/(p-1)\ZB\bigr) 
\cross (\ZB/2\ZB)}$ 
 ([25], VI, Satz 5.1).
\medskip
\noindent {\sl  Question.} Does there exist a commutative  field  -- or even a ring -- $K$ such that
${\rm Br}_{\rm Gal}(K)$ is strictly contained in ${\rm Br}(K)$?
\medskip
This problem looks rather difficult to solve since very few is known in general
about Brauer groups.  A natural candidate seems to be an
infinite field $K$ of characteristic 2. Indeed, there exists then 
a central simple  $K$-algebra that is
not Galois over $K$ (namely the quaternion algebra 
${\displaystyle {H_{a,b} =  [{{a,b}\over K})} }$ over $K$ of counter-example 1.2). 
The technics involved here however
do not allow to show that under suitable conditions none of the representatives in
the class $[H_{a,b}] \in {\rm Br}(K)$ is Galois.

\medskip
\noindent 3.5. {\sl The category ${\goth Gal}_R(G)$ and base change}. 
\smallskip
\noindent
Denote by $\pi _{alg}^{gr}$ the functor from
${\goth Gal}$ to ${\goth Gpf}^{\rm o}\cross {\goth Alg}_k$
defined by $\pi _{alg}^{gr}(G,\psi : R \lr S) = (G,R)$. 
Introduce the fibre 
category, denoted by
${\goth Gal}_R(G)$, of $\pi _{alg}^{gr}$ over $(G,R)$.
Its objects are the $G$-Galois extensions  $\psi : R \lr S$, with
$R$ and $G$ fixed. A
morphism from $ \psi : R \lr S$ to $\psi  : R \lr S'$
in ${\goth Gal}_R(G)$
is a 
$G$-equivariant morphism of $R$-rings
$\varphi : S \lr S'$.
The full subcategories   
${\goth Galstr}_R(G)$, ${\goth Galcent}_R(G)$, ${\goth GalAzum}_R(G)$
and ${\goth Galcom}_R(G)$
of ${\goth Gal}_R(G)$ are defined in an obvious way. 
\medskip

 Le Bruyn, 
van den Bergh and van Oystaeyen proved that if $\psi : R \lr S$ is
a Galois extension of {\sl commutative} rings for some finite group $G$, and if 
 $T$ is an arbitrary $R$-algebra, then $T \lr S\ot _RT$ is also a $G$-Galois
extension ([20], Lemma II.5.1.13). Next we state another version
of base change, assuming that the initial Galois extension 
is centralizing or central and the algebra $T$ is commutative.

\medskip
\noindent {\bf  Proposition 3.6.} {\sl  (Base change over an extension of commutative rings). 
Let $R$ be a commutative algebra and 
$G$ a finite group.
A morphism of commutative algebras
$\varrho : R\lr T$ induces two {\sl ``extension of the scalars by $T$"}  functors
$\varrho _! : {\goth Galcent}_R(G) \lr {\goth Galcent}_T(G)$ and
$ _!\varrho : {\goth Galcent}_R(G) \lr {\goth Galcent}_T(G) $
(resp.
$\varrho _! : {\goth GalAzum}_R(G) \lr {\goth GalAzum}_T(G)$ and 
$_!\varrho : {\goth GalAzum}_R(G) \lr {\goth GalAzum}_T(G) $)
defined by $$\varrho _! (\psi : R \lr S) = (\psi _T : T \lr T\ot _RS)
\ \ \ \ {\hbox {\sl and}} \ \ \ \ _!\varrho (\psi : R \lr S) = (_T\psi  : T \lr S\ot _RT),$$
where, for any $t\in T$,
 $$\psi _T(t) = t\ot 1 \ \ \ \ {\hbox {\sl and}} \ \ \ \ _T\psi (t) = 1\ot t.$$}

\noindent {\sl Proof.}  We show the proposition for instance for the left 
extension of the scalars by $T$. 
The proof paraphrases  the one of
Theorem 3.1. 
The group $G$ acts on $U = T\ot _RS$ by
$g(t\ot s)
= t\ot g(s),$ for $g \in G$ and $t \ot s \in U.$
The abelian group $U = T\ot _RS$ is clearly a left $T$-module. 
Let it also be a right $T$-module by the formula $(t\ot s)\cdot t' = tt' \ot s$,
for $t, t' \in T$ and $s\in S$. Since $T$ is commutative, $U$ becomes  symmetric
as a $T$-module. 

 1) It is clear that $T \cong T \ot _RR \subseteq U^G$. A descent argument analogue 
to the one used in
Theorem 3.1 shows that if $u\in U^G$,
then $u \in {\rm Im}(\id _T\ot \tr _{S/R})$ = $T\ot _RR \cong T, $ 
because $u = (\id _T \ot \tr _{S/R})\bigl((1\ot s_0)u\bigr)$, for any fixed  $s_0\in S$
which verifies  $\tr _{S/R}(s_0) = 1$. Indeed, set ${\displaystyle u = \sum_{i = 1}^mt_i \ot s_i}$.
For any $g \in G$, one has ${\displaystyle u = \sum_{i = 1}^mt_i \ot g(s_i)}$. So
$$\eqalign{ u &= (1\ot 1)\cdot u = (\id _T \ot \tr _{S/R})\bigl((1\ot s_0)u\bigr) = 
(\id _T \ot \tr _{S/R})\bigl(\sum_{i = 1}^m(t_i\ot s_0s_i)\bigr) \cr & =
\sum_{g \in G}\sum_{i = 1}^m\bigl(t_i\ot g(s_0) g(s_i)\bigr) 
 = \sum_{g \in G}\bigl(1 \ot g(s_0)\bigr)\bigl(\sum_{i = 1}^mt_i\ot g(s_i)\bigr) 
= \sum_{g \in G}\bigl(1 \ot g(s_0)\bigr)u.\cr}$$

 2) The morphism $\Gamma : U\ot _TU \lr U(G)$ defined by
$$\Gamma (u\ot u') =
\sum_{g\in G}u g(u')\delta_g$$
is an isomorphism of left $T$-modules as a composition 
 $\beta \circ (\id_T \ot \Gamma _{S}) \circ \alpha$, where
$$\alpha :U\ot _RU \lr T\ot _RS\ot _RS \quad {\rm and} 
\quad \beta : T\ot _R\bigl(S(G)\bigr)\lr U(G)$$ are the canonical isomorphisms.

 3) At last, $T$ is a direct summand in $U$,
since $R$ is a direct summand in $S$ and the base change functor
$T\ot _R-$ preserves direct limits, hence direct sums. 
 \dm

\medskip
\noindent {\bf  Applications 3.7.} 

\item{--} {\sl Localization.}
Let $\psi : R\lr S$ be a centralizing $G$-Galois extension 
and $\Sigma$ be a multiplicative subset of $R$. Then $G$ acts on 
the localized algebra $\Sigma ^{-1}S = \Sigma ^{-1}R\ot _RS$ of $S$ by $\Sigma$ via the formula
$$g({s \over u}) = {g(s) \over u},$$
for $s \in S$ and $u \in \Sigma .$ The localized morphism $\Sigma ^{-1}\psi : \Sigma ^{-1}R \lr \Sigma ^{-1}S$
remains a centralizing $G$-Galois extension. In particular, for any prime ideal $\wp \in {\rm Spec}(R)$
(res\-pectively maximal ideal ${\goth m} \in {\rm Max}(R)$),
the morphism of local rings $\psi _\wp : R_\wp \lr S_\wp$ (respectively $\psi _{\goth m} : R_{\goth m} \lr S_{\goth m}$) is $G$-Galois.

\smallskip

\item{--} {\sl  The case of commutative rings.} 
The opposite functor ${\displaystyle \pi _{alg}^{opp} : {\goth Galcomm}^{opp} \lr {\goth Aff}_k}$ 
to the functor
$\pi _{alg}$ defined on $ {\goth Galcomm}$ with values in ${\goth Algcom}_k$ 
enables us to build over any affine $k$-scheme $\Spec R $ the fibre
category  of $\pi _{alg}$ denoted ${\displaystyle {\goth Galcomm}^{opp}_R}$. 
Every morphism of affine $k$-schemes
 $f : \Spec T \lr \Spec R$ (or
equivalently, every morphism of commutative algebras $\theta : R \lr T$)
induces then a base change functor 
${\displaystyle f^* = \theta_!^{opp}: {\goth Galcomm}^{opp}_R \lr {\goth Galcomm}^{opp}_T}$.
Observe that the collection of these data is coherent in the following sense:
if $f$ and $g$ are composable morphisms of affine schemes, 
there exists a natural transformation 
$\phi _{f,g} : (f\circ g)^* \simeq g^*\circ f^*.$

\medskip

\noindent {\bf  Corollary 3.8.} {\sl The assignment 
$R \longmapsto {\rm Br}_{\rm Gal}(R)$ defines a functor ${\rm Br}_{\rm Gal}$
from the category ${\goth Algcom}_k$ of commutative algebras to the category of groups.
For each morphism
$\varrho : R\lr T$ in ${\goth Algcom}_k$, the map ${\rm Br}_{\rm Gal}(\varrho ) $ is the
homomorphism of groups ${\rm Br}_{\rm Gal}(R) \lr$ $ {\rm Br}_{\rm Gal}(T)$
induced by the extension of the scalars $\varrho : R \lr T$.}

\medskip

\Dem Let $\varrho : R\lr T$ be a morphism in ${\goth Algcom}_k$. The assignment $A \longmapsto A \ot _RT$ induces a 
$\ZB$-homomorphism ${\rm Br}(\varrho ) : {\rm Br}(R) \lr {\rm Br}(T)$ ([9], Theorem 8.9). 
By base change over the commutative algebra 
morphism $R \lr T$ (Proposition 3.6), 
the Brauer functor ${\rm Br} : {\goth Algcom}_k \lr {\goth Gp}$ 
restricts  to ${\rm Br}_{\rm Gal}(\varrho ) : {\rm Br}_{\rm Gal}(R) \lr {\rm Br}_{\rm Gal}(T)$.
\dm

\medskip 
\noindent {\sl  Example 3.9.}
 Let $K$ be a global field and $\Omega$ be a complete set of valuations of 
$K$. For each  $v \in \Omega$, let $K_v$ be the completion
of $K$ at $v$. The collection 
of group homomorphisms ${\rm Br}_{\rm Gal}(K) \lr {\rm Br}_{\rm Gal}(K_v)$ deduced from $K \lr K_v$
provides  a map $$\displaystyle{\beta : {\rm Br}_{\rm Gal}(K) \lr  \bigoplus_{v \in \Omega}{\rm Br}_{\rm Gal}(K_v)}$$
which is a monomorphism of groups.
Indeed, the injectivity of $\beta$ immediately results from the short exact sequence 
$$0 \lr {\rm Br}(K) \lr \bigoplus_{v \in \Omega}{\rm Br}(K_v) \lr \QB/\ZB \lr 0$$
obtained in the classical theory of 
Hasse-Brauer-Noether-Albert ([11], Chap. 13).

\medskip 
\noindent 3.6. {\sl  The relative Brauer-Galois group.}
\smallskip
\noindent
Let $\varrho : R\lr T$  be a morphism
of commutative algebras and let ${\rm Br}(\varrho ) $ be the
homomorphism of groups ${\rm Br}(R) \lr {\rm Br}(T)$
induced by $\varrho$.
The kernel of ${\rm Br}(\varrho )$ is called the 
{\sl relative Brauer group of $\varrho$}
and is denoted by ${\rm Br}(T/R)$. If $A$ is an $R$-Azumaya algebra such that its Brauer class $[A]$
lies in  ${\rm Br}(T/R)$, one says that $T$ {\sl splits} $A$.
Corollary 3.8 allows us to define the {\sl relative Brauer-Galois group of $\varrho : R \lr T$}, 
denoted by ${\rm Br}_{\rm Gal}(T/R)$,
to be the kernel of ${\rm Br}_{\rm Gal}(\varrho )$; it is a subgroup of  ${\rm Br}(T/R)$.
\medskip 
\noindent {\sl  Remarks.}

\noindent 1.-- When $R = K$ is a commutative field and $A$ is a $K$-Azumaya algebra, 
$A$ is split by some finite commutative Galois extension
of $K$, say $L$, with Galois group $G = {\rm Gal}(L/K)$ ([28], Corollary 3.15). This strong result
does not allow any conclusion about conditions whether $K \lr A$ should be  Galois  or not.
 
\noindent 2.-- If $\psi : R\lr S$ is a finite Galois extension of commutative rings with Galois group $G =
{\rm Gal}(S/R)$ 
such that
the Picard group ${\rm Pic}(S)$ of $S$ is trivial, then there  is a well known
cohomological characterization of the relative Brauer  group ([28], Theorem 7.12):
the  group
${\rm Br}(S/R)$
is isomorphic to the Galois $2$-cohomology group $H^2(G, S^{\lcross})$.
Hence ${\rm Br}_{\rm Gal}(S/R)$ corresponds to a subgroup $H^2_{{\rm Gal}}(G, S^{\lcross})$
contained in  $H^2(G, S^{\lcross})$.

\vskip 20pt
\noindent {\bf Acknowledgments:} The author wishes to thank Lars Kadison for 
many useful comments on an early version of this paper.

\medskip
\medskip
\medskip
\vfill
\eject

\noindent{\bf REFERENCES}

\medskip

\item{[1]} M. A{\ninecmr USLANDER} -- O. G{\ninecmr OLDMAN}:
The Brauer group of a commutative ring, 
{\sl Trans. Amer. Math. Soc.} {\bf 97}, 
({\oldstyle 1960}), 367  --  409.
\medskip

\item{[2]} A. B{\ninecmr LANCHARD}: 
{\sl Les corps non commutatifs},
Presses Universitaires de France, Paris ({\oldstyle 1972}).
\medskip

\item{[3]} N. B{\ninecmr OURBAKI}: {\sl Alg\`ebre commutative}, Chap. 1 et 2, 
Hermann, Paris, ({\oldstyle 1961}).
\medskip

\item{[4]} S. C{\ninecmr AENEPEEL}  --  G. M{\ninecmr ILITARU}  --   S. Z{\ninecmr HU}:  
{\sl Frobenius and separable functors for generalized module categories and nonlinear equations},
Lecture Notes in Mathematics 1787, Springer-Verlag, Berlin --
Heidelberg -- New York -- Barcelona -- Hong-Kong -- London -- Milan --
Paris -- Tokyo ({\oldstyle 2002}).
\medskip

\item{[5]} S. U. C{\ninecmr HASE} -- D. K. H{\ninecmr ARRISON} --
A. R{\ninecmr OSENBERG}: Galois theory and cohomology of 
commutative rings,
{\sl Memoirs of the A.M.S.,}
Number 52  ({\oldstyle 1965}).
\medskip

\item{[6]} S. U. C{\ninecmr HASE} -- M. E. S{\ninecmr WEEDLER}:
{\sl Hopf algebras and Galois theory}, Lecture
Notes in Math. 97, Springer-Verlag, Berlin --
Heidelberg -- New York ({\oldstyle 1969}). 
\medskip

\item{[7]} F. D{\ninecmr E}M{\ninecmr EYER}  --  E. I{\ninecmr NGRAHAM}:  
{\sl Separable algebras over commutative rings},
Lecture Notes in Mathematics 181, Springer-Verlag, Berlin --
Heidelberg -- New York ({\oldstyle 1971}).
\medskip

\item{[8]} P. K. D{\ninecmr RAXL}:  
{\sl Skew fields},
London Mathematical Society Lecture Note Series {\bf 81}, Cambridge University Press,
Cambridge
({\oldstyle 1983}).
\medskip

\item{[9]} B. F{\ninecmr ARB} --  R.K. D{\ninecmr ENNIS}:
{\sl Noncommutative algebra},
Graduate Texts in Mathematics, Vol. 144, 
Springer-Verlag, Berlin -- Heidelberg -- New York ({\oldstyle 1993}).
\medskip

\item{[10]} C. G{\ninecmr REITHER}:  
{\sl Cyclic Galois extensions of commutative rings},
Lecture Notes in Math. 1534, 
Springer-Verlag, Berlin -- Heidelberg  ({\oldstyle 1992}).
\medskip

\item{[11]} A. J. H{\ninecmr AHN}:  
{\sl Quadratic algebras, Clifford algebras, and arithmetic Witt groups},
Universitext, 
Springer-Verlag, Berlin --  Heidelberg -- New-York ({\oldstyle 1993}).
\medskip

\item{[12]} L. K{\ninecmr ADISON}:
{\sl New examples of Frobenius extensions}, University Lecture Series, Vol. 14,
American Mathematical Society, Providence R.I. 
({\oldstyle 1999}).
\medskip

\item{[13]} T. K{\ninecmr ANZAKI}:
On commutator rings and Galois theory of separable algebras, 
{\sl Osaka J. Math.} {\bf 1} 
({\oldstyle 1964}), 103  --  115; Correction, {\sl ibid.} {\bf 1} 
({\oldstyle 1964}), 253.
\medskip

\item{[14]} C. K{\ninecmr ASSEL}:
{\sl Quantum groups},
Graduate Texts in Mathematics, Vol. 155, 
Springer-Verlag, Berlin -- Heidelberg -- New York ({\oldstyle 1995}).
\medskip

\item{[15]} M. A. K{\ninecmr NUS}:  
{\sl Quadratic and hermitian forms over rings},
 Grundlehren der mathema\-tischen Wissenschaften 294, 
Springer-Verlag, Berlin --
Heidelberg -- New York ({\oldstyle 1991}).
\medskip

\item{[16]} M. A. K{\ninecmr NUS}  --  A. M{\ninecmr ERKURJEV}
 --  M. R{\ninecmr OST}  --  J.-P. T{\ninecmr IGNOL}:  
{\sl The Book of Involutions}, American Mathematical Society, Colloquium
Publications, Volume 44, Providence ({\oldstyle 1998}).
\medskip

\item{[17]} M. A. K{\ninecmr NUS}  --  M. O{\ninecmr JANGUREN}:  
{\sl Th\'eorie de la descente et alg\`ebres d'Azumaya}, 
Lecture Notes in Mathematics 389, Springer-Verlag, Berlin --
Heidelberg -- New York ({\oldstyle 1974}).
\medskip

\item{[18]} H. F. K{\ninecmr REIMER} -- M. T{\ninecmr AKEUCHI}: Hopf algebras and 
Galois extensions of an algebra,
{\sl Indiana Univ. Math. J.}, Vol. {\bf 30} ({\oldstyle 1981}), 
675  --  692.
\medskip

\item{[19]} S. L{\ninecmr ANG}: {\sl Algebra}, Third edition,
Addison Wesley ({\oldstyle 1997}).
\medskip

\item{[20]} L. {\ninecmr LE} B{\ninecmr RUYN}  --  M. V{\ninecmr AN DEN}
B{\ninecmr ERGH}  --  F. V{\ninecmr AN} O{\ninecmr YSTAEYEN}:
{\sl Graded orders}, Birkh\"auser, Boston -- Basel ({\oldstyle 1988}).
\medskip

\item{[21]} B. A. M{\ninecmr AGURN}: {\sl An algebraic introduction to $K$-Theory}, 
Encyclopedia of Mathematics and its Applications, Cambridge University Press,
Cambridge ({\oldstyle 2002}).
\medskip

\item{[22]} {\cyr A. S. M{\ninecyr ERKUR{\char 94}EV}:
O gomomorfizme normennogo vyqeta stepeni dva, 
{\cyri Dokl. AN SSSR}, t.} {\bf 264}, {\cyr {\char 125}} 3
({\oldstyle 1981}), 542  --  547. {\sl English translation:} 
A. S. M{\ninecmr ERKUR'EV}: 
On the norm residue symbol of degree 2, 
{\sl Soviet Math. Dokl.}, Vol. {\bf 24(3)}
({\oldstyle 1981}), 546  --  551.
\medskip

\item{[23]} {\cyr A. S. M{\ninecyr ERKUR{\char 94}EV}  --  A. A. S{\ninecyr USLIN}:
$K$-kogomologii mnogoobrazi{\char 26} Se\-ve\-ri-Brau{\char 11}ra i gomomorfizm normennogo vyqeta, 
{\cyri Izv. AN SSSR, Ser. Mat.}, t.} {\bf 46}, {\cyr {\char 125}} 5
({\oldstyle 1983}), 1011  --  1046, 1135  --  1136. {\sl English translation:} 
A. S. M{\ninecmr ERKUR'EV}  --   A. A. S{\ninecmr USLIN}: 
$K$-cohomology of Severi-Brauer varieties and the norm residue homomorphism, 
{\sl Math. USSR -- Izv.}, Vol. {\bf 21(2)}
({\oldstyle 1983}), 307  --  340.
\medskip

\item{[24]} J. M{\ninecmr ILNOR}: {\sl Introduction to algebraic $K$-Theory}, Annals of Mathematics Studies,
Princeton University Press, Princeton ({\oldstyle 1971}).
\medskip

\item{[25]} J. N{\ninecmr EUKIRCH}:  
{\sl Algebraische Zahlentheorie},
Springer-Verlag, Berlin -- Heidelberg -- New York -- London ({\oldstyle 1992}).
\medskip

\item{[26]} P. N{\ninecmr USS}:
Noncommutative descent and nonabelian cohomology, 
{\sl $K$-Theory} {\bf 12} 
({\oldstyle 1997}), 23  --  74.
\medskip

\item{[27]} P. N{\ninecmr USS}:
Extensions galoisiennes non commutatives~: normalit\'e,
cohomologie non ab\'elienne, 
{\sl Comm. Algebra} {\bf 28 (7)} 
({\oldstyle 2000}), 3223  --  3251.
\medskip

\item{[28]} M. O{\ninecmr RZECH}  --  C. S{\ninecmr MALL}:  
{\sl The Brauer group of commutative rings},
Lecture Notes in pure and applied Mathematics 11, Marcel Dekker Inc.,  New York ({\oldstyle 1975}).
\medskip

\item{[29]} R. S. P{\ninecmr IERCE}:
{\sl Associative algebras},
Graduate Texts in Mathematics, Vol. 88, 
Springer-Verlag, Berlin -- Heidelberg -- New York ({\oldstyle 1982}).
\medskip

\item{[30]} J. R{\ninecmr OSENBERG}:  
{\sl Algebraic $K$-theory and its applications}, 
Graduate Texts in Mathematics, Vol. 147, 
Springer-Verlag, Berlin --
Heidelberg -- New York ({\oldstyle 1994}).
\medskip

\item{[31]} J.-P. S{\ninecmr ERRE}:  
{\sl  Corps locaux}, Troisi\`eme \'edition corrig\'ee,
Hermann, Paris  ({\oldstyle 1968}).

\bye